\documentclass[conference]{IEEEtran}
\IEEEoverridecommandlockouts
\usepackage{amsmath,amsfonts,amsthm} % Apparently amssymb conflicts with the next one
\usepackage[cmintegrals]{newtxmath} % This is the recommended style.
\usepackage{algorithmic}
\usepackage{graphicx}
\graphicspath{{images/}}
\usepackage{textcomp}
\usepackage{xcolor}
\def\BibTeX{{\rm B\kern-.05em{\sc i\kern-.025em b}\kern-.08em
    T\kern-.1667em\lower.7ex\hbox{E}\kern-.125emX}}
    
\usepackage{biblatex}
\usepackage{multirow} % For colorbar which spans two rows
\usepackage{array}
    
\newcommand{\ie}[0]{\textit{i.e.}}
\newcommand{\eg}[0]{\textit{e.g.}}  

\newtheorem{theorem}{Theorem}[section]

\theoremstyle{definition}
\newtheorem{definition}[theorem]{Definition}

\theoremstyle{remark}
\newtheorem{remark}[theorem]{Remark}

\usepackage{tikz}
\usetikzlibrary{shapes.geometric, arrows, positioning, calc}
\tikzstyle{input} = [rectangle, rounded corners, text centered, draw=black, fill=yellow!20, minimum width = 2cm, minimum height=0.8cm,]
\tikzstyle{treeA} = [rectangle, text centered, draw=black, fill=blue!30] % Node for tree A
\tikzstyle{treeB} = [rectangle, text centered, draw=black, fill=red!30] % Node for tree B
\tikzstyle{arrow} = [thick,->,>=stealth] % Connecting 'arrow'
\tikzstyle{title} = [text width=1cm, text centered] % Titles for directions

\newcommand{\xyplane}[4]{% color, opacity, x-coordinate, y-coordinate
\filldraw[thick,color=black,fill=#1,opacity=#2] (0,0,0) -- (#3,0,0) -- (#3,0,#4) -- (0,0,#4) -- cycle;
}
\newcommand{\yzplane}[4]{% color, opacity, y-coordinate, z-coordinate
\filldraw[thick,color=black,fill=#1,opacity=#2] (0,0,0) -- (#3,0,0) -- (#3,#4,0) -- (0,#4,0) -- cycle;
}
\newcommand{\xzplane}[4]{% color, opacity, x-coordinate, z-coordinate
\filldraw[thick,color=black,fill=#1,opacity=#2] (0,0,0) -- (0,#3,0) -- (0,#3,#4) -- (0,0,#4) -- cycle;
}

\DeclareMathOperator{\R}{\mathbb{R}}
\DeclareMathOperator{\N}{\mathbb{N}}
\DeclareMathOperator{\Z}{\mathbb{Z}}
\DeclareMathOperator{\C}{\mathbb{C}}

\newcommand{\Wre}{\psi_{\text{Re}}} % Real part of the wavelet function
\newcommand{\Wim}{\psi_{\text{Im}}} % Imaginary part of wavelet function
\newcommand{\Sre}{\varphi_{\text{Re}}} % Real part of the scaling function
\newcommand{\Sim}{\varphi_{\text{Im}}} % Imaginary part of scaling function
\newcommand{\Gre}{\gamma_{\text{Re}}} % Real part of placeholder function gamma
\newcommand{\Gim}{\gamma_{\text{Im}}} % Imaginary part of placeholder function gamma

\renewcommand{\H}{\mathcal{H}} % Hilbert transform.
\newcommand{\I}{\mathbf{i}} % Imaginary unit function for different orthants

\renewcommand{\Re}[1]{\text{Re}\left(#1\right)} % Real part
\renewcommand{\Im}[1]{\text{Im}\left(#1\right)} % Imag. part

\renewcommand{\vec}[1]{\boldsymbol{#1}} % vectors in boldface
\newcommand{\f}{\vec{f}}
\newcommand{\m}{\vec{m}}
\newcommand{\db}{\vec{d}}
\newcommand{\vi}{\vec{v}}
\newcommand{\etab}{\vec{\eta}}

\DeclareMathOperator{\Radon}{\mathcal{R}}
\DeclareMathOperator{\RadonD}{\vec{\mathcal{R}}}
\DeclareMathOperator{\ComplexWTD}{\vec{\mathcal{C}}}
\DeclareMathOperator{\WaveTD}{\vec{\mathcal{W}}}
\newcommand{\proj}[1]{\text{proj}_{#1}} % Projection operator

\begin{document}

\title{4D Dual-Tree Complex Wavelets for Time-Dependent Data}

\author{\IEEEauthorblockN{Tatiana A.~Bubba}
\IEEEauthorblockA{\textit{Dept. of Mathematics and Statistics} \\
\textit{University of Helsinki}\\
Helsinki, Finland \\
tatiana.bubba@helsinki.fi}
\and
\IEEEauthorblockN{Tommi Heikkil\"{a}}
\IEEEauthorblockA{\textit{Dept. of Mathematics and Statistics} \\
\textit{University of Helsinki}\\
Helsinki, Finland \\
tommi.heikkila@helsinki.fi}
\and
\IEEEauthorblockN{Samuli Siltanen}
\IEEEauthorblockA{\textit{Dept. of Mathematics and Statistics} \\
\textit{University of Helsinki}\\
Helsinki, Finland \\
samuli.siltanen@helsinki.fi}
\thanks{TAB acknowledges support by the Academy of Finland postdoctoral grant, decision number 330522. TH and SS acknowledge support by the Academy of Finland Project 310822.
All authors acknowledge partial support by Academy of Finland through the Finnish Centre of Excellence in Inverse Modelling and Imaging 2018–2025, decision number 312339.}
}

\maketitle

\begin{abstract}
The dual-tree complex wavelet transform (DT-$\C$WT) is extended to the 4D setting. Key properties of 4D DT-$\C$WT, such as directional sensitivity and shift-invariance, are discussed and illustrated in a tomographic application. The inverse problem of reconstructing a dynamic three-dimensional target from X-ray projection measurements can be formulated as 4D space-time tomography. The results suggest that 4D DT-$\C$WT offers simple implementations combined with useful theoretical properties for tomographic reconstruction. 
\end{abstract}


\section{Introduction}

\IEEEPARstart{A}{} wide selection of multiscale methods have been introduced in recent years for representing and processing  multidimensional signals. It is well-known that classical wavelets~\cite{Mallat2009} are not optimal for anisotropic data in dimensions two and higher. On the other hand, they offer simple and relatively fast implementations (especially considering the curse of dimensionality), strong theoretical properties (bases and orthogonality) and thorough theoretical understanding.

Complex-valued wavelets, and in particular the dual-tree implementation originally introduced by N.~Kingsbury~\cite{kingsbury1998dual} and extended to 3D in~\cite{chen2011efficient}, utilize most of these advantages. Additionally, they provide directional sensitivity and shift-invariance with a simpler construction than those of curvelets~\cite{candes2002new} or shearlets~\cite{kutyniok2012shearlets}.

These nice features also ease the extension of the dual-tree complex wavelets to higher dimensions, especially 4D, where concepts like specific directions and even visualization are obviously difficult. In some sense the natural world is 4-dimensional (3 spatial dimensions and time) and, more concretely, a wide variety of different and interesting 4D data arises from spectral imaging, geospatial applications,  computer graphics, and more.

This motivates the main contribution of our work, namely the extension of the construction of the dual-tree complex wavelet transform (DT-$\C$WT) to 4D. Our Matlab implementation of the 4D dual-tree complex wavelet transform and its inverse called, respectively, \texttt{dualtree4} and \texttt{idualtree4} (after the 2D and 3D implementations of similar names) is available on GitHub~\cite{heikkila2021dualtree4}. 

We demonstrate the feasibility of the DT-$\C$WT for 4D applications by applying it to the inverse problem of reconstructing a changing volume over time from a collection of X-ray images. In this 3D+time dynamic computed tomography (CT), the 4D DT-$\C$WT helps to overcome the ill-posedness of the inverse problem via \textit{regularization} \cite{Engl1996}. The above-mentioned favorable theoretical properties  allow details and edges be preserved over time in the reconstructions, even when the measurements are very sparsely collected (only 30 projection views). The 4D DT-$\C$WT also outperform real-valued wavelet transform computationally in this application.

While 4D real-valued wavelets have been considered in applications (\eg{}, \cite{ansari2020robust,jahanian20084d}), to our knowledge a 4D complex-valued wavelet system has not been proposed before. In addition different extensions of wavelets using quaternions~\cite{chan2004quaternion} (which are 4D in a different sense) and hypercomplex numbers~\cite{chan2004directional} have been introduced but the actual implementations have so far been limited to 2D and 3D setting.

The rest of this paper is organized as follows. In section~\ref{sec:4DDTCWT} we introduce the 4D DT-$\C$WT, after briefly revising the construction of the DT-$\C$WT in 2D. Properties of the (4D) DT-$\C$WT, like shift-invariance and directional sensitivity, are shortly illustrated in section~\ref{sec:properties}. In section~\ref{sec:applications} we apply the 4D DT-$\C$WT as a regularizer to the ill-posed problem of 4D dynamic CT
problem: we test our model on both a simulated and a physical phantom. Finally, we draw some conclusions in section~\ref{sec:conclusions}.

\section{Implementation and algorithm}
\label{sec:4DDTCWT}
The name dual-tree comes from the original implementation for 1D signals where the two real-valued discrete wavelet transform (DWT) trees are used side by side to obtain the real and imaginary parts of the complex wavelet coefficients for all scales of the decomposition. In 2D (and higher dimensions) the two DWTs are no longer as separated due to the way higher dimensional wavelets are usually constructed but its 1D components still share this original design. Moreover, a dual construction can still be used for a simple and efficient implementation of complex-valued wavelets in higher dimensions, including 4D. 
We first cover this method in 2D, then move to 4D and finally consider both the inverse and adjoint of the transform. We also formalize all this by defining the associated operators.

%%%%%%%%%%
% 2D & 3D
%%%%%%%%%%
\subsection{Constructing DT-$\C$WT in 2D and 3D}
\label{subsec:CWT2D}
Before going into the details of the 4D dual-tree complex wavelet transform (DT-$\C$WT) we begin by briefly discussing the construction of the 2D version. In some sense the key changes happen when the complex wavelets are extended from one dimension to two and from then on it is simply a matter of accounting the larger number of filters and their permutations. For a more detailed account on DT-$\C$WT, we refer to the papers by the original authors \cite{kingsbury1998dual,kingsbury2001complex} (which include the 2D transform) and their work on extending it to 3D \cite{chen2011efficient}.

Similarly to the 2D (real-valued) discrete wavelet transform~\cite{Mallat2009}, to define the 2D DT-$\C$WT 
we can use any 1D (complex-valued) mother wavelet $\psi (x)= \psi_{\C} (x) = \Wre(x) + i\Wim(x)$ (associated with a high-pass filter $H$) and scaling function $\varphi(x) = \varphi_{\C} (x)= \Sre(x) + i \Sim (x)$ (associated with a low-pass filter $L$). By taking their tensor product and switching their role along the directions $x$ and $y$, we obtain the 2D scaling and wavelet functions. For example, the 2D wavelet whose both directions use the wavelet function (typically denoted by $HH$) is given by:
\begin{align}
    \psi_{\C}(x,y) &:= \psi(x) \psi(y) \nonumber \\
    &= [\Wre(x) + i\Wim(x)] \times [\Wre(y) + i\Wim(y)] \\
    &= \Wre(x) \Wre(y) - \Wim(x) \Wim(y) \nonumber \\
    &\quad + i\left(\Wre(x) \Wim(y) + \Wim(x) \Wre(y) \right). \label{eq:2Dorthant1}
\end{align}
If $\Wre$ and $\Wim$ form (approximately) a Hilbert transform pair
$\Wim = \H(\Wre)$ (\ie{}, they are $90^{\circ}$ out of phase with each other) then $\psi$ is (approximately) analytic and $\widehat{\psi}(\hat{x})$ vanishes for $\hat{x} < 0$. This means that the 2D wavelet $\psi_{\C}(x,y)$ is only supported on the positive orthant of the frequency domain ($\hat{x}, \hat{y} > 0$). Let's denote it by $\psi_{O1}(x,y)$ and define a second wavelet:
\begin{align}
\psi_{O2}(x,y) &:= \overline{\psi(x)} \psi(y) \nonumber \\
&= \Wre(x) \Wre(y) + \Wim(x) \Wim(y) \nonumber \\
    &+ i\left(\Wre(x) \Wim(y) - \Wim(x) \Wre(y) \right), \label{eq:2Dorthant2}
\end{align}
where $\overline{\psi(x)} = \Wre(x) - i\Wim(x)$ is the complex conjugate of $\psi(x)$. Then $\psi_{O2}$ is supported on the second orthant of the frequency domain ($\hat{x} < 0, \hat{y} > 0$). Similarly, we could define wavelets for the other two orthants (where $\hat{y} < 0$) but if we only wish to apply our wavelet transform to real valued functions $f$, then it is not necessary since:
\begin{equation}
    \langle f, \psi_{O3} \rangle = \langle f, \psi(x) \overline{\psi(y)} \rangle = \overline{\langle \overline{f}, \overline{\psi(x)} \psi(y) \rangle} = \overline{\langle f, \psi_{O2} \rangle},
\end{equation}
meaning the coefficients corresponding to $O3$ are complex conjugates of $O2$.
The symmetry is analogous between $\psi_{O1}$ and $\psi_{O4}$. This shows that we only need complex conjugation in the $x$-component.

Next, we still need to consider the other complex wavelet configurations given by
\begin{align*}
    \psi(x) \varphi(y) & \qquad (HL \text{ wavelet}), \\
    \varphi(x) \psi(y) & \qquad (LH \text{ wavelet}), \\
    \varphi(x) \varphi(y) & \qquad (LL \text{ wavelet}).
\end{align*}
The last one is known as the 2D scaling function or father wavelet, which we will consider separately later. In total, we have three 2D complex wavelet configurations and, for each one, we also need to consider the complex conjugate on the $x$-component which means there are $3 \cdot 2 = 6$ different 2D complex wavelet functions. However, it is not necessary to consider all 6 explicitly as we did with the $HH$ wavelet. Instead, we will introduce some new notation which will be particularly useful later on in the 4D setting.

First, note that both pairs of 1D functions $\lbrace \Wre, \Sre \rbrace$ and $\lbrace \Wim, \Sim \rbrace$ constitute regular 1D real-valued wavelet systems which are only connected through the Hilbert transform. Therefore, they are associated with two different pairs of low-pass, high-pass filters: $H_a$ and $L_a$ for the real part; $H_b$ and $L_b$ for the imaginary one. The construction of these so-called ``q-shift'' wavelet filters is thoroughly explained and motivated in \cite{kingsbury2001complex,kingsbury2003design}, to which we refer the reader for a more detailed discussion. 

In the original 1D DT-$\C$WT the two wavelet systems would produce the two independent halves of a tree-like structure, whence the name: tree $a$ would only produce the real part while tree $b$ only the imaginary one.
This is no longer true in 2D as we see, for example, in equation~\eqref{eq:2Dorthant1}. The real part is a sum of two functions, one of which is purely from tree $a$ and the second one is purely from tree $b$. The imaginary part is a sum of two terms made by mixing both trees. How the filters from the two trees are multiplied and added is the key to a simpler implementation.
By comparing equations \eqref{eq:2Dorthant1} and \eqref{eq:2Dorthant2} we notice that both 2D wavelets are computed by summing up the same 4 terms and only the signs change due to the imaginary unit changing from positive to negative. In fact, the orthants alone determine the signs no matter which filter is used.

Therefore, we introduce the following notation. Define the real-valued terms $P_{\iota}$, where the multi-index $\iota \in \lbrace aa, ab, ba, bb \rbrace$ denotes from which tree the filters along the $x$ and $y$-directions are chosen from. To be precise these terms should be unique to each wavelet configuration: for example, for the $HL$ wavelet the precise notation should be $P_{H_a L_a}, P_{H_a L_b}, P_{H_b L_a}$ and $P_{H_b L_b}$. However, if we ease the notation by dropping the explicit dependence on the filters, we can generalize computations by considering only the dependence on the tree. The actual filter types can be inferred from the generated wavelet.
With this notation, \textit{all} wavelets (\ie{}, $LH$, $HL$ and $HH$) in the two orthants with $\hat{y}>0$ have the following form:
\begin{align*}
    \psi_{O1}(x,y) &= P_{aa} - P_{bb} + i \left( P_{ab} + P_{ba} \right), \\
    \psi_{O2}(x,y) &= P_{aa} + P_{bb} + i \left( P_{ab} - P_{ba} \right)
\end{align*}
and to compute them, we only need the corresponding $P_{\iota}$ for all $\iota \in \lbrace aa, ab, ba, bb \rbrace$. 
Since computing each $P_{\iota}$ term amounts to the same complexity as computing wavelet coefficients with any real-valued 2D wavelet transform, this provides a considerable simplification for the implementation of 2D DT-$\C$WT.

The same approach is also used to calculate each $P_\iota$ corresponding to the scaling function ($LL$), but the different orthants are not considered explicitly since the complex-valued coefficients are not extracted from the $P_\iota$ terms.
Instead, these are stored in an alternating pattern where values of $P_{aa}$ are stored on even columns and rows, $P_{ab}$ on even columns and odd rows and so on.
This produces a single larger set of real-valued scaling coefficients which is then passed on as the input for the next decomposition level.

At the final decomposition level, the complex-valued scaling coefficients for the two orthants are computed just like the wavelet coefficients, that is:
\begin{align*}
    \varphi_{O1}(x,y) &= P_{aa} - P_{bb} + i \left( P_{ab} + P_{ba} \right), \\
    \varphi_{O2}(x,y) &= P_{aa} + P_{bb} + i \left( P_{ab} - P_{ba} \right).
\end{align*}

Finally, to invert the 2D DT-$\C$WT, we reverse the operations above and obtain each term $P_{\iota}$ from the respective complex-valued coefficients:
\begin{align*}
    P_{aa} &= \frac{1}{2} \Re{\psi_{O1} + \psi_{O2}}, \qquad
    P_{bb} = \frac{1}{2} \Re{\psi_{O1} - \psi_{O2}}, \\
    P_{ab} &= \frac{1}{2} \Im{\psi_{O1} + \psi_{O2}}, \qquad
    P_{ba} = \frac{1}{2} \Im{\psi_{O1} - \psi_{O2}}.
\end{align*}
Hence, also the reconstruction can be carried out as in any 2D real-valued wavelet system.

Constructing the DT-$\C$WT system in 3D follows similarly. We will not go into details here since the construction and an application of the system are thoroughly discussed in~\cite{chen2011efficient}. Just to build intuition to then generalize to 4D, the wavelet and scaling functions are constructed as a tensor product of 3 complex-valued components for the $x$, $y$ and $z$-directions respectively. In 3D we have $2^3 = 8$ configurations, 1 for the 3D scaling function ($LLL$) and 7 for the different wavelet functions. Similarly to the 2D case, complex conjugated components are needed in two directions to cover the negative parts of the respective frequency domain while in the third direction the covering is obtained by symmetry, if only real-valued inputs are used.
In total $7 \cdot 2^{3-1} = 28$ different 3D wavelets at each decomposition level are needed. 

%%%%%%%%%
% 4D
%%%%%%%%%
\subsection{Constructing DT-$\C$WT in 4D}
To construct the 4D DT-$\C$WT, we extend the approach described in section~\ref{subsec:CWT2D} to four dimensions, where the different directions are denoted by $x$, $y$, $z$ and $t$. Also in 4D the different wavelet configurations have a separable construction, using tensor products of 1D (dual-tree) complex-valued wavelet or scaling function for any given direction. This yields $2^4 = 16$ configurations, 1 of which corresponds to the 4D scaling function while the other 15 are wavelets.

Again, the 1D complex-valued wavelet and scaling functions have, by construction, frequency support only on one half of the domain. Thus, in order to obtain complete frequency tiling, a complex conjugated function needs to be included for each dimension as well. Once more, with real-valued inputs one order of symmetry is obtained in the Fourier domain and hence, for one dimension, conjugated functions are not needed: let us fix this to be the fourth dimension, corresponding to $t$. 

As a result, the total number of 4D wavelet functions used is $15 \cdot 2^{(4-1)} = 15 \cdot 8 = 120$ for each decomposition level. Luckily, 120 unique filters are not explicitly needed since we can use the same trick seen in the 2D setting to obtain the 8 directional orthants.

All mother wavelet functions have the same form:
\begin{align*}
    \psi_{O\zeta}(x,y,z,t) :&= [\Gre(x) + \I_x(\zeta) \Gim(x)] \\
    & \times [\Gre(y) + \I_y(\zeta) \Gim(y)] \\
    & \times [\Gre(z) + \I_z(\zeta) \Gim(z)] \\
    & \times [\Gre(t) + \I_t(\zeta) \Gim(t)]
\end{align*}
where, depending on the desired wavelet, the functions $\Gre$ and $\Gim$ can be wavelet ($\Wre, \Wim$) or scaling functions ($\Sre, \Sim$) in any of the 15 configurations. Here, $\I_d(\zeta) = \pm i$ determines the sign of the imaginary unit $i$ for each dimension $d$ based on the orthant $\zeta$. Thanks to the symmetry with real valued inputs, we can fix $\I_t(\zeta) = +i$ for any $\zeta$. All values of $\I_d$ are listed in table~\ref{tab:orthants}. 

\begin{figure}
    \centering
    \begin{tikzpicture}
    \node (origin) at (0,0) {}; % Empty center point for easier placement
    % input box
    \node (input) [input, left=0.02cm of origin] {\large 4D filter};
    \draw[thick] (input) -- (0,0);
    
    % x-direction
    % a
    \node (a) [treeA, above=4.0cm of origin] {Tree $a$};
    % b
    \node (b) [treeB, below=4.0cm of origin] {Tree $b$};
    \draw [thick,<->,>=stealth] (a) -- node[anchor=west,pos=0.9] {$\I_x$} (b);
    
    \node (xtitle) [title, above=0.2cm of a] {$x$ ``rows''}; % x-direction title
    
    % y-direction
    % aa
    \node (aa) [treeA, above right=1.7cm and 0.3cm of a] {Tree $a$};
    \draw [arrow] (a) -| (aa);
    % ab
    \node (ab) [treeB, below right=1.7cm and 0.3cm of a] {Tree $b$};
    \draw [arrow] (a) -| node[anchor=west,pos=0.8] {$\I_y$} (ab);
    % ba
    \node (ba) [treeA, above right=1.7cm and 0.3cm of b] {Tree $a$};
    \draw [arrow] (b) -| node[anchor=west,pos=0.8] {$\I_x$} (ba);
    % bb
    \node (bb) [treeB, below right=1.7cm and 0.3cm of b] {Tree $b$};
    \draw [arrow] (b) -| node[anchor=west,pos=0.8] {$\I_x \I_y$} (bb);
    
    \node (ytitle) [title, above=0.2cm of aa] {$y$ ``columns''}; % y-direction title
    
    % z-direction
    % aaa
    \node (aaa) [treeA, above right=0.6cm and 0.3cm of aa] {Tree $a$};
    \draw [arrow] (aa) -| (aaa);
    % aab
    \node (aab) [treeB, below right=0.6cm and 0.3cm of aa] {Tree $b$};
    \draw [arrow] (aa) -| node[anchor=west,pos=0.8] {$\I_z$} (aab);
    % aba
    \node (aba) [treeA, above right=0.6cm and 0.3cm of ab] {Tree $a$};
    \draw [arrow] (ab) -| node[anchor=west,pos=0.8] {$\I_y$} (aba);
    % abb
    \node (abb) [treeB, below right=0.6cm and 0.3cm of ab] {Tree $b$};
    \draw [arrow] (ab) -| node[anchor=west,pos=0.8] {$\I_y \I_z$} (abb);
    % baa
    \node (baa) [treeA, above right=0.6cm and 0.3cm of ba] {Tree $a$};
    \draw [arrow] (ba) -| node[anchor=west,pos=0.8] {$\I_x$} (baa);
    % bab
    \node (bab) [treeB, below right=0.6cm and 0.3cm of ba] {Tree $b$};
    \draw [arrow] (ba) -| node[anchor=west,pos=0.8] {$\I_x \I_z$} (bab);
    % bba
    \node (bba) [treeA, above right=0.6cm and 0.3cm of bb] {Tree $a$};
    \draw [arrow] (bb) -| node[anchor=west,pos=0.8] {$\I_x$} (bba);
    % bbb
    \node (bbb) [treeB, below right=0.6cm and 0.3cm of bb] {Tree $b$};
    \draw [arrow] (bb) -| node[anchor=west,pos=0.8] {$\I_x \I_z$} (bbb);
    
    \node (ztitle) [title, above=0.2cm of aaa] {$z$ ``slices''}; % z-direction title
    
    % t-direction
    % aaaa
    \node (aaaa) [treeA, above right=0.25cm and 0.3cm of aaa] {Tree $a$};
    \draw [arrow] (aaa) -| (aaaa);
    \node[right=0.1cm of aaaa] {$= P_{aaaa}$};
    % aaab
    \node (aaab) [treeB, draw=purple,very thick, dashed, below right=0.25cm and 0.3cm of aaa] {Tree $b$};
    \draw [arrow] (aaa) -| node[anchor=west,pos=0.8] {$\I_t$} (aaab);
    \node[right=0.1cm of aaab] {$= P_{aaab}$};
    % aaba
    \node (aaba) [treeA, draw=purple,very thick, dashed, above right=0.25cm and 0.3cm of aab] {Tree $a$};
    \draw [arrow] (aab) -| node[anchor=west,pos=0.8] {$\I_z$} (aaba);
    \node[right=0.1cm of aaba] {$= P_{aaba}$};
    % aabb
    \node (aabb) [treeB, below right=0.25cm and 0.3cm of aab] {Tree $b$};
    \draw [arrow] (aab) -| node[anchor=west,pos=0.8] {$\I_z \I_t$} (aabb);
    \node[right=0.1cm of aabb] {$= P_{aabb}$};
    % abaa
    \node (abaa) [treeA, draw=purple,very thick, dashed, above right=0.25cm and 0.3cm of aba] {Tree $a$};
    \draw [arrow] (aba) -| node[anchor=west,pos=0.8] {$\I_y$} (abaa);
    \node[right=0.1cm of abaa] {$= P_{abaa}$};
    % abab
    \node (abab) [treeB, below right=0.25cm and 0.3cm of aba] {Tree $b$};
    \draw [arrow] (aba) -| node[anchor=west,pos=0.8] {$\I_y \I_t$} (abab);
    \node[right=0.1cm of abab] {$= P_{abab}$};
    % abba
    \node (abba) [treeA, above right=0.25cm and 0.3cm of abb] {Tree $a$};
    \draw [arrow] (abb) -| node[anchor=west,pos=0.8] {$\I_y \I_z$} (abba);
    \node[right=0.1cm of abba] {$= P_{abba}$};
    % abbb
    \node (abbb) [treeB, draw=purple,very thick, dashed, below right=0.25cm and 0.3cm of abb] {Tree $b$};
    \draw [arrow] (abb) -| node[anchor=west,pos=0.8] {$\I_y \I_z \I_t$} (abbb);
    \node[right=0.1cm of abbb] {$= P_{abbb}$};
    
    % more t-direction
    % baaa
    \node (baaa) [treeA, draw=purple,very thick, dashed, above right=0.25cm and 0.3cm of baa] {Tree $a$};
    \draw [arrow] (baa) -| node[anchor=west,pos=0.8] {$\I_x$} (baaa);
    \node[right=0.1cm of baaa] {$= P_{baaa}$};
    % baab
    \node (baab) [treeB, below right=0.25cm and 0.3cm of baa] {Tree $b$};
    \draw [arrow] (baa) -| node[anchor=west,pos=0.8] {$\I_x \I_t$} (baab);
    \node[right=0.1cm of baab] {$= P_{baab}$};
    % baba
    \node (baba) [treeA, above right=0.25cm and 0.3cm of bab] {Tree $a$};
    \draw [arrow] (bab) -| node[anchor=west,pos=0.8] {$\I_x \I_z$} (baba);
    \node[right=0.1cm of baba] {$= P_{baba}$};
    % babb
    \node (babb) [treeB, draw=purple,very thick, dashed, below right=0.25cm and 0.3cm of bab] {Tree $b$};
    \draw [arrow] (bab) -| node[anchor=west,pos=0.8] {$\I_x \I_z \I_t$} (babb);
    \node[right=0.1cm of babb] {$= P_{babb}$};
    % bbaa
    \node (bbaa) [treeA, above right=0.25cm and 0.3cm of bba] {Tree $a$};
    \draw [arrow] (bba) -| node[anchor=west,pos=0.8] {$\I_x \I_y$} (bbaa);
    \node[right=0.1cm of bbaa] {$= P_{bbaa}$};
    % bbab
    \node (bbab) [treeB, draw=purple,very thick, dashed, below right=0.25cm and 0.3cm of bba] {Tree $b$};
    \draw [arrow] (bba) -| node[anchor=west,pos=0.8] {$\I_x \I_y \I_t$} (bbab);
    \node[right=0.1cm of bbab] {$= P_{bbab}$};
    % bbba
    \node (bbba) [treeA, draw=purple,very thick, dashed, above right=0.25cm and 0.3cm of bbb] {Tree $a$};
    \draw [arrow] (bbb) -| node[anchor=west,pos=0.8] {$\I_x \I_y \I_z$} (bbba);
    \node[right=0.1cm of bbba] {$= P_{bbba}$};
    % bbbb
    \node (bbbb) [treeB, below right=0.25cm and 0.3cm of bbb] {Tree $b$};
    \draw [arrow] (bbb) -| node[anchor=west,pos=0.8] {$\I_x \I_y \I_z \I_t$} (bbbb);
    \node[right=0.1cm of bbbb] {$= P_{bbbb}$};
    
    \node (ttitle) [title, above=0.2cm of aaaa] {$t$ ``time step''}; % t-direction title
    
    \end{tikzpicture}
    \caption{Illustration of the tree-like structure which determines the terms $P_{\iota}$, with $\iota \in \lbrace aaaa, aaab, \ldots, bbbb \rbrace$, of every wavelet filter and the combined imaginary units coming from each dimension. Cells highlighted with dashed boundary form the imaginary part of the final output.}
    \label{fig:filterTree}
\end{figure}
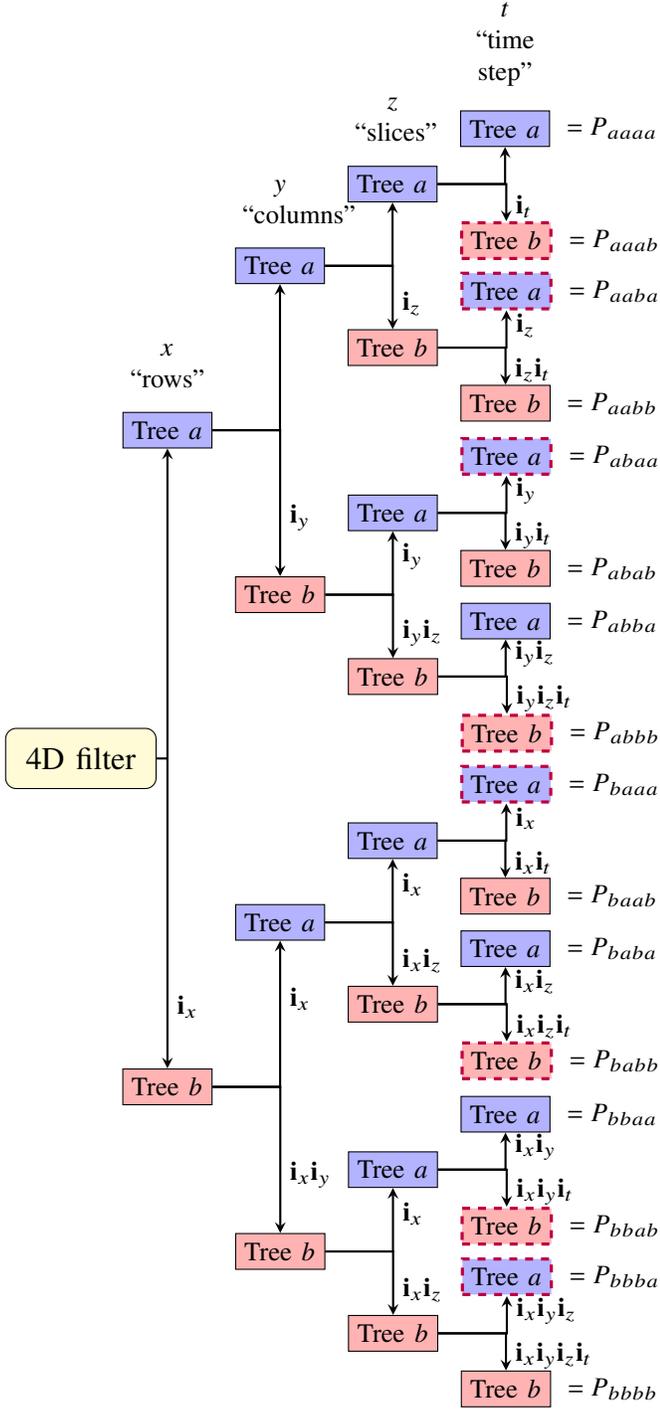

Calculating this product always gives a sum of 16 terms, 8 for the real part (\ie{}, even number of $\I_d$'s) and 8 for the imaginary part (\ie{}, odd number of $\I_d$'s). The sign of each term is determined solely by the product of the $\I_d$'s. It is clear that by simply changing the signs of the 16 terms all 8 orthants can be covered by any given wavelet. Choosing how the imaginary units affects each term can be easily seen from the tree-like structure in Figure~\ref{fig:filterTree}: each time a filter is chosen from the tree $b$ for direction $d$, each subsequent term on that branch is multiplied by $\I_d$.

Since an even number of ``imaginary'' wavelets produces the real part of the coefficient, neither wavelet tree is purely imaginary (or real) valued. For this reason we denote them by $a$ and $b$ instead. Then, for tree $a$ we have high-pass filter $H_a$ and low-pass filter $L_a$ corresponding to wavelet and scaling functions $\Wre$ and $\Sre$, respectively. Similarly for tree $b$ we have high-pass filter $H_b$ and low-pass filter $L_b$ corresponding to $\Wim$ and $\Sim$, respectively. Individually they produce orthogonal real-valued wavelet systems and they are only connected by the Hilbert transform pairing $\Wim = \H(\Wre)$.

\begin{figure}
    \centering
    \begin{tikzpicture}[scale=2.5]
    \node (origin) at (0,0,0) {}; % Center point
    % Due to obvious overlap, there is a specific order of drawing the regions
    % Behind yz plane
    \xyplane{gray}{0.3}{-1}{-1} % xy: back left
    \xzplane{gray}{0.3}{-1}{-1} % xz: back bottom
    \xyplane{gray}{0.3}{1}{-1} % xy: back right
    \xzplane{gray}{0.3}{1}{-1} % xz: back top
    % Orthant #: y z x
    \node at (.8,.8,-.9) {$O2$};
    \node at (-.8,.8,-.9) {$O4$};
    \node at (.8,-.8,-.9) {$O6$};
    %\node at (-.8,-.9,-.8) {$O8$};
    \draw[thick,->] (-1.2,-.6,.3) -- (-1.01,-.6,0) node[pos=0,left] {$O8$};
    % yz plane
    \yzplane{gray}{0.3}{1}{1} % yz: top left
    \yzplane{gray}{0.3}{-1}{1} % yz: top right
    \yzplane{gray}{0.3}{1}{-1} % yz: bottom right
    \yzplane{gray}{0.3}{-1}{-1} % yz: bottom left
    % In front of yz plane
    \xzplane{gray}{0.5}{-1}{1} % xz: front bottom
    \xyplane{gray}{0.5}{-1}{1} % xy: front left
    \xzplane{gray}{0.5}{1}{1} % xz: front top
    \xyplane{gray}{0.5}{1}{1} % xy: front right
    % Orthant #: y z x
    \node at (.6,.7,.4) {$O1$};
    \node at (-.6,.7,.4) {$O3$};
    \node at (.6,-.6,.4) {$O5$};
    \node at (-.6,-.6,.4) {$O7$};
    
    % Coordinate arrows
    \draw[very thick,->,>=stealth] (origin) -- (0,1.2,0) node[above] {$\hat{z}$};
    \draw[very thick,->,>=stealth] (origin) -- (1.2,0,0) node[above] {$\hat{y}$};
    \draw[very thick,->,>=stealth] (origin) -- (0,0,1.3) node[left] {$\hat{x}$};
    \end{tikzpicture}
    \caption{Illustration of the different orthants in the  3D $(\hat{x},\hat{y},\hat{z})$-Fourier space; the $\hat{t}$-dimension is left out for clarity.}
    \label{fig:orthants}
\end{figure}
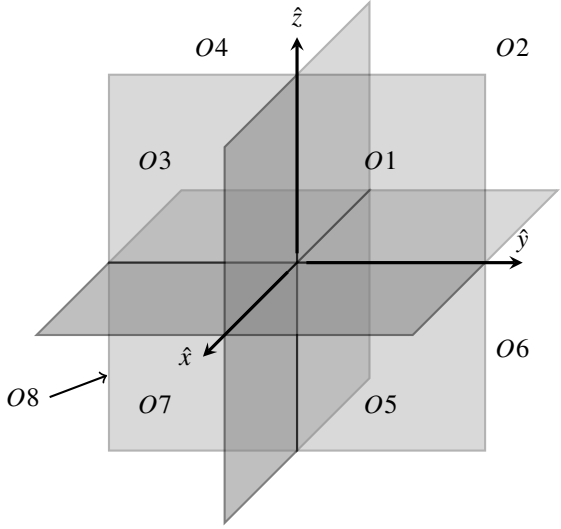

As an example consider the complex-valued wavelet denoted by $HHHH$ and obtained using a high-pass filter $H$ in every direction.
To uniquely identify this wavelet, we need to further differentiate whether the filter $H$ is from tree $a$ or $b$. As in 2D, let us denote these terms by $P_{\iota}$ where the multi-index $\iota$ marks the tree for each of the four dimensions, namely $\iota \in \lbrace aaaa, aaab, \ldots, bbbb \rbrace$. This is precisely the structure illustrated in Figure~\ref{fig:filterTree}. Then, each $HHHH$ wavelet is given by 16 terms, all of which are computed using the different high-pass filters from the two trees.
The signs for the imaginary units in the different orthants of the frequency domain are given in Table~\ref{tab:orthants}. See also Figure~\ref{fig:orthants} for an illustration of the orthants in the $(\hat{x},\hat{y},\hat{z})$-Fourier space; the $\hat{t}$-dimension is left out for clarity.

\begin{table}
    \centering
    \caption{Imaginary unit values for each of the orthants. Negative imaginary unit corresponds to complex conjugation of the complex wavelet component in that direction.}
    \label{tab:orthants}
    \begin{tabular}{c|cccc}
        Orthant & $\I_x$ & $\I_y$ & $\I_z$ & $\I_t$ \\ \hline
        $O1$ & $+i$ & $+i$ & $+i$ & $+i$ \\
        $O2$ & $-i$ & $+i$ & $+i$ & $+i$ \\
        $O3$ & $+i$ & $-i$ & $+i$ & $+i$ \\
        $O4$ & $-i$ & $-i$ & $+i$ & $+i$ \\
        $O5$ & $+i$ & $+i$ & $-i$ & $+i$ \\
        $O6$ & $-i$ & $+i$ & $-i$ & $+i$ \\
        $O7$ & $+i$ & $-i$ & $-i$ & $+i$ \\
        $O8$ & $-i$ & $-i$ & $-i$ & $+i$ \\
    \end{tabular}
\end{table}

For the real and imaginary parts of any wavelet in, \eg{}, the first orthant this yields the following expressions:
\begin{align}
\begin{split}
    \Re{\psi_{O1}} &= \frac{1}{2} \big(P_{aaaa} - P_{aabb} - P_{abab} - P_{abba} \\
    & - P_{baab} - P_{baba} - P_{bbaa} + P_{bbbb}\big), \\[0.5em]
    \Im{\psi_{O1}} &= \frac{1}{2} \big(P_{aaab} + P_{aaba} + P_{abaa} - P_{abbb} \\
    & + P_{baaa} - P_{babb} - P_{bbab} - P_{bbba}\big),
\end{split}    
\label{eq:psiO1}
\end{align}
whereas the real and imaginary parts of any wavelet in the second orthant (where $\I_x = -i$) are be given by:
\begin{align}
\begin{split}
    \Re{\psi_{O2}} &= \frac{1}{2} \big(P_{aaaa} - P_{aabb} - P_{abab} - P_{abba} \\
    & + P_{baab} + P_{baba} + P_{bbaa} - P_{bbbb}\big), \\[0.5em]
    \Im{\psi_{O2}} &= \frac{1}{2} \big(P_{aaab} + P_{aaba} + P_{abaa} - P_{abbb} \\
    & - P_{baaa} + P_{babb} + P_{bbab} + P_{bbba}\big).
\end{split}    
\label{eq:psiO2}
\end{align}
As we can see the terms coming from the bottom half of the tree in Figure~\ref{fig:filterTree} have their signs changed because, by construction, the imaginary unit $\I_x$ is always present in those terms. All terms are also multiplied by a factor of $\frac{1}{2}$ to lower the frame bound of the wavelet system. This is not strictly necessary and some other options for the normalization are mentioned in subsection~\ref{ssec:inverse}.

To obtain all the other possible configurations, we simply reiterate the same procedure. Namely, for every configuration of the low-pass and high-pass filters we obtain 16 terms denoted by $P_{\iota}$ where $\iota$ keeps track on which tree each filter was chosen from. By changing the signs of these terms we can obtain the real and imaginary parts of any ``directional orthant''. Therefore, the sign of each $P_{\iota}$ and how its values are computed are independent of each other.

\begin{remark}[first decomposition level]
\label{rm:FirstLevDec}
The first decomposition level differences from the others in that it uses just one low-pass and high-pass filter which correspond to a biorthogonal wavelet system. Here, instead of tree~$a$ and tree~$b$ the final output consists of the odd and even values of the filters convolved with the odd and even values of the input. This method is simpler than using two sets of q-shift filters and faster to compute thanks to shorter filters. This is especially advantageous in higher dimensions where most of the computational load is on the first decomposition level. For example in 4D every subsequent level is only $\frac{1}{16}th$ of the size of the previous one.

However, using just one set of filters does not work properly beyond the first level due to the different sampling rate, therefore the q-shift filters are required. Nevertheless, once the convolutions produce the $P_{\iota}$ terms, the rest of the computations are carried out identically in every decomposition level.

Furthermore, one could also consider simply discarding these first level details coefficients, gaining faster computations at the cost of an imperfect final reconstruction, due to the partial missing information encoded by these detail coefficients. This option is offered both by Matlab's built-in 3D DT-$\C$WT~\cite{chen2011efficient} and also our \texttt{dualtree4} implementation~\cite{heikkila2021dualtree4}. 
\end{remark}

\medskip

We conclude the subsection by formally introducing the definition of complex wavelet transform $\mathcal{C}$ which acts as \emph{analysis  operator}, \ie{}, it maps (decomposes) any input $f$ to its complex wavelet coefficients.

\begin{definition}
Let $\psi^{(\kappa,\zeta)}(x,y,z,t), \ \kappa = 0,1,\ldots,15, \ \zeta = 1,\ldots,8$ denote the different mother wavelets (based on the configurations and orthants), including for $k=0$ the scaling function $\varphi = \psi^{(0)}$, with a slight abuse of notation. Let $f \in L^2(\R^4)$ be a real-valued function.
Then $\mathcal{C}$ is defined to be the linear mapping:
\begin{equation}
\label{eq:analysis}
    \mathcal{C}: f(x,y,z,t) \, \longmapsto \, 
    c(j,m;\kappa,\zeta) = \langle f, \psi_{j,m}^{(\kappa,\zeta)} \rangle
\end{equation} 
where $\psi_{j,m}^{(\kappa,\zeta)} = 2^{-2j} \psi^{(\kappa,\zeta)}(2^{-j}\cdot - m )$ with $(j,m) \in \N_0 \times \Z^4$.
\end{definition}

Notice that, since in the numerical setting the scale is in practice limited $j \leqslant J$ and we must also include translates of the scaling function (namely, for $\kappa = 0$ with $\varphi^{(\zeta)} = \varphi^{(0, \zeta)}$). This can be done, for example, by defining  
\begin{align*}
    c(J,m;\kappa,\zeta) = \langle f, \varphi^{(\zeta)}_{J,m} \rangle = \left\langle f, 2^{-2J} \varphi^{(\zeta)}(2^{-J}\cdot -m) \right\rangle
\end{align*}
for the different orthants $\zeta = 1,\ldots,8$.
%\txtp{Quick question: do we need $k$ here? Because now we have $\kappa =1,\ldots, 15$ and this covers only the mother wavelets, not the father right?}\txtb{I think denoting the father wavelet as "16th wavelet" is the easiest option.}
As usual, the maximum decomposition level is bounded by the resolution of the data: $2^J \leqslant \min\lbrace N_x, N_y, N_z, N_t \rbrace$. For practical reasons \texttt{dualtree4} only works when each $N_x, N_y, N_z$ and $N_t$ is even.

\subsection{Inverting 4D DT-$\C$WT} 
\label{ssec:inverse}
Inverting the dual-tree complex wavelet decomposition is a very straight forward process once the terms $P_{\iota}$ are separated from the complex valued coefficients. Let $\Re{\psi_{O\zeta}}$ and $\Im{\psi_{O\zeta}}$ denote, respectively, the real and imaginary parts of some complex wavelet at scale $j$ and orthant $\zeta$. Similarly to the 2D case, we can compute the corresponding $P_{\iota}$ as follows:
\begin{align}
\begin{split}
    P_{aaaa} = \frac{1}{4}&\Big(\text{Re}\big(\psi_{O1} + \psi_{O2} + \psi_{O3} + \psi_{O4} \\
    &+ \psi_{O5} + \psi_{O6} + \psi_{O7} + \psi_{O8} \big)\Big), \\
    P_{aabb} = \frac{1}{4}&\Big(\text{Re}\big(\psi_{O1} - \psi_{O2} - \psi_{O3} - \psi_{O4} \\
    &+ \psi_{O5} + \psi_{O6} + \psi_{O7} + \psi_{O8} \big)\Big),
\end{split}  
\label{eq:InvTrasf1}
\end{align}
and so on for the terms which were summed for the real part in~\eqref{eq:psiO1}. Similarly the terms which were summed for the imaginary part are given by 
\begin{align}
\begin{split}
    P_{aaab} = \frac{1}{4}&\Big(\text{Im}\big(\psi_{O1} + \psi_{O2} + \psi_{O3} + \psi_{O4} \\
    &+ \psi_{O5} + \psi_{O6} + \psi_{O7} + \psi_{O8} \big)\Big), \\
    P_{abbb} = \frac{1}{4}&\Big(\text{Im}\big(- \psi_{O1} - \psi_{O2} + \psi_{O3} + \psi_{O4} \\
    &+ \psi_{O5} + \psi_{O6} - \psi_{O7} - \psi_{O8} \big)\Big),
\end{split}  
\label{eq:InvTrasf2}    
\end{align}
and analogously for the remaining $P_{\iota}$ terms.
The division by $4$ is required since in the decomposition step each term $P_{\iota}$ is divided by $2$ and here we obtain $8 \cdot \frac{1}{2}P_{\iota} = 4P_{\iota}$ for the desired term while the rest cancel out. Another option would be to use a uniform normalization of $\frac{1}{\sqrt{8}}$ for both the decomposition and the reconstruction steps which would produce a Parseval frame but also be computationally slightly more expensive than multiplying by a fraction.

From this point onward the reconstruction is carried out just like with any DWT. For levels $j > 1$ the reconstruction filters $\widetilde{H}_a, \widetilde{L}_a, \widetilde{H}_b$ and $\widetilde{L}_b$ are ``time-reversed'' (\ie{}, the 1D filters are mirrored) versions of the respective decomposition filters. For $j=1$ the reconstruction filters $\widetilde{H}, \widetilde{L}$ are the associated dual filters of the biorthogonal wavelet system.

\medskip

We end the subsection by formally defining the inverse 4D DT-$\C$WT $\mathcal{C}^{-1}$ which allows to reconstruct the original signal from its DT-$\C$WT coefficients.

\begin{definition}
Let $\kappa = 0, \ldots ,15, \, \zeta = 1,\ldots,8$ and $(j,m) \in \N_0 \times \Z^4$.
The inverse complex wavelet transform $\mathcal{C}^{-1}$ is given by
\begin{align}
    \mathcal{C}^{-1}&: c(j,m;\kappa,\zeta) \, \longmapsto \, f(x,y,z,t) \nonumber \\
    f &= \sum_{j} \sum_m \sum_{\kappa, \zeta} c(j,m;\kappa, \zeta)\widetilde{\psi}^{(\kappa,\zeta)}_{j,m}.
\end{align}
\end{definition}
Here, $\widetilde{\psi}$ marks the dual wavelet function of the biorthogonal wavelet system used at $j=1$ as mentioned in remark~\ref{rm:FirstLevDec}. For $j \geqslant 2$ these are the same as for the analysis operator. In the numerical setting the scaling function is once again included with $\kappa = 0$.

\subsection{Adjoint 4D DT-$\C$WT}

In some applications (such as the one we propose in section~\ref{sec:applications}) the adjoint $\mathcal{C}^*$ of the complex wavelet transform $\mathcal{C}$ is required in place of the inverse. The adjoint of the analysis operator is also known as the \emph{synthesis operator}. Since the orthogonal wavelet systems used for levels $j \geqslant 2$ use the same filters (just time-reversed) for the inverse, the adjoint is the inverse but scaled by $\frac{1}{2}$ (instead of the normalization factor $\frac{1}{4}$ in equations~\eqref{eq:InvTrasf1} and~\eqref{eq:InvTrasf2}): namely, it has the same normalization factor of decomposition operator. However, for the first level the dual filters $\widetilde{H}, \widetilde{L}$ also need to be replaced by the time-reversed decomposition filters $H, L$.

This produces a fairly accurate approximation of the adjoint operator and in our implementation is available by using the parameter ``adjoint'' when calling the function \texttt{idualtree4}~\cite{heikkila2021dualtree4}. Further improvement could be obtained by a more detailed consideration of the boundary conditions of the discrete convolution, as explained in~\cite{folberth2016efficient}, but we leave this to future work.

It is worth mentioning that since this particular implementation does not constitute a Parseval frame but a tight frame with frame bound $u = 2$, this bound is also present in the adjoint. Hence, the largest eigenvalue of the normal operator $\mathcal{CC}^*$ is $2^2$.

\section{Properties}
\label{sec:properties}
Since the dual-tree complex wavelet system is constructed using two real-valued DWT systems side by side, computationally it is at least $2^4 = 16$ times as demanding as using real-valued discrete wavelet transform of similar filter lengths. However, the dual-tree complex wavelets exhibit many appealing properties (lacking in the real-valued DWT) which make them a tempting option many tasks.

\subsection{Shift-invariance}
While real-valued wavelets are well suited for many applications, their implementation is in general sensitive to small translations in the input. This means that the DWT coefficients from data which have been slightly shifted can significantly differ from those of the non-shifted data.

This is not the case for DT-$\C$WT. Since the real and imaginary parts of the dual-tree complex wavelet are in quadrature (\ie{}, $90^\circ$ difference in phase) and the absolute value of the wavelet is not oscillatory, errors caused by shifts are in general less severe. In fact, aiming for shift-invariant wavelets leads precisely to complex-valued wavelets: shift-invariance can be numerically confirmed using various filters, as shown in~\cite{kingsbury2001complex}. % Confusing paper with numerical examples 
As an example, figure~\ref{fig:timecone} demonstrates how shift-invariance in DT-$\C$WT, coupled with its directional sensitivity (see subsection~\ref{subsec:directionality}), helps preserving edges over time. 

For a particular class of complex-valued wavelets, called modulated wavelets, it is possible to formally prove that the errors caused by shifts are optimally small~\cite{barri2012near}. We leave the extension of this result to the 4D dual-tree complex wavelets presented in this paper to future work.

\subsection{Directionality}
\label{subsec:directionality}
One of the main drawbacks of real-valued DWT is the lack of ability of capturing directional information in 2-dimensions and beyond. This was the main reason for introducing multidimensional systems like curvelets~\cite{candes2002new} or shearlets~\cite{kutyniok2012shearlets}. From a theoretical perspective, complex-valued wavelets share certain limitations of real-valued wavelet systems\footnote{For example, the asymptotic decay rate remains $\mathcal{O}(N^{-1})$ in 2D~\cite{kutyniok2012shearlets} and $\mathcal{O}(N^{-\frac{1}{d-1}})$ in $d$-dimensions for $d-1$ dimensional edges, which is known to be suboptimal in terms of best nonlinear $N$-term approximation.}, given that the scaling is still isotropic and there is no explicit encoding of directionality. However, in practice, it can be seen that dual-tree complex wavelets can capture directional information across a fixed number of orientations per scale.

\begin{figure}[t]
    \centering
    \includegraphics[width=0.24\textwidth]{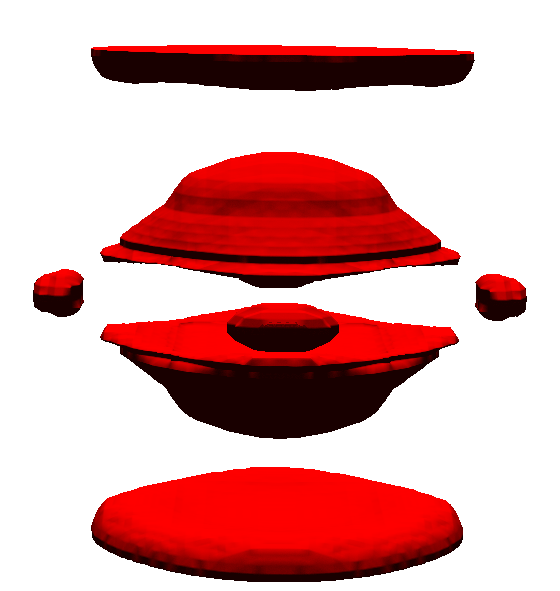}
    \includegraphics[width=0.24\textwidth]{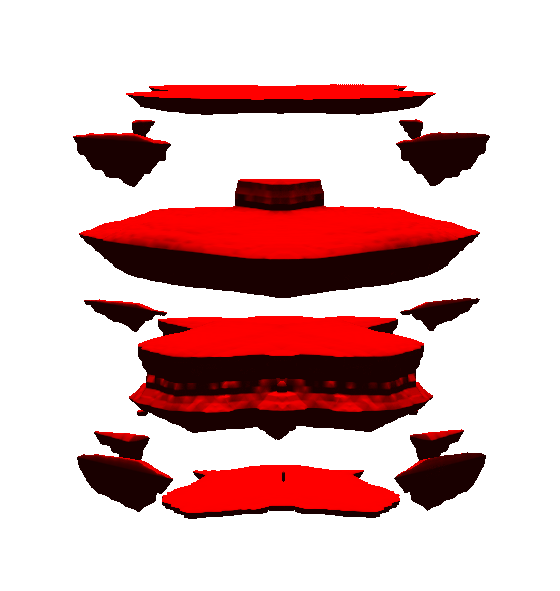}
    \caption{3D isosurface rendering of a ball growing over time reconstructed using only the coefficients corresponding to complex-valued (left) or real-valued (right) wavelets aligned vertically (along the $z$-axis). Only the middle time step ($t=16$) is shown here.}
    \label{fig:LLHLrecon}
\end{figure}

Indeed, with dual-tree complex wavelets, details in different parts of the spatial domain are analyzed by wavelets supported in different orthants of the Fourier domain. This ``one-sided frequency support''
results in a major selectivity (compared to DWT) in representing singularities which eventually entails the ability to naturally encode some directionality.

In figure~\ref{fig:LLHLrecon}, we demonstrate this by comparing the reconstruction of a simple 3D ball growing over time using the coarsest scale ``LLHL''-wavelet coefficients of both dual-tree complex wavelets and (Daubechies 2) real-valued wavelets.  

It is clear from figure~\ref{fig:LLHLrecon} that dual-tree complex wavelets (left) produce a remarkably cleaner representation of the edges of the ball in the vertical direction and the reconstruction remains symmetric. 

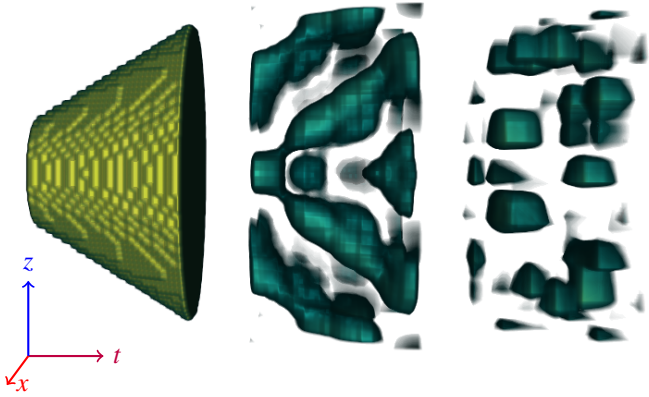
\begin{figure}[t]
    \centering
    \begin{tikzpicture}
        \node[inner sep=0pt] (timecone) at (-0.5,-0.1)
        {\includegraphics[width=0.48\textwidth]{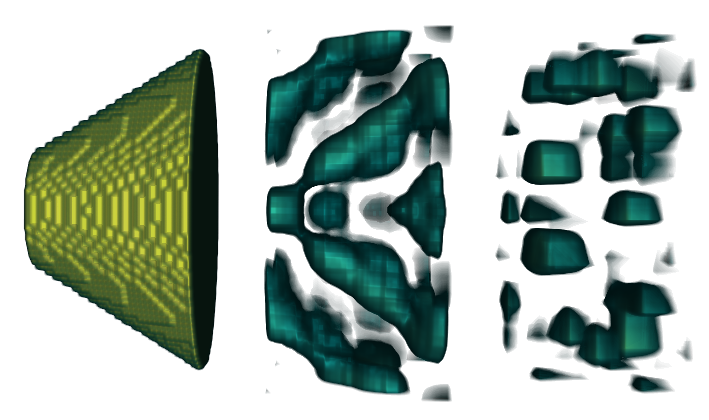}}; % Image
        % Coordinate axes
        \coordinate (O) at (-5,-3,-1.2); % Move this point to move all axes
        \draw[->,purple,thick] (O) -- ($(O) + (1,0,0)$) node[anchor=west,pos=1] {$t$};
        \draw[->,blue,thick] (O) -- ($(O) + (0,1,0)$) node[anchor=south,pos=1] {$z$};
        \draw[->,red,thick] (O) -- ($(O) + (0.1,0,1)$) node[anchor=west,pos=1] {$x$};
    \end{tikzpicture}
    \caption{3D ``time-cone'' rendering of a central $xz$-slice of a ball growing over time (left) reconstructed using only the coefficients which correspond to a complex-valued (middle) or real-valued (right) wavelets aligned diagonally (along the $zt$-plane).}
    \label{fig:timecone}
\end{figure}

Instead, real-valued wavelets (right) result in a reconstruction with jagged and unintuitive edges and the overall rounded shape of the growing ball seems to be lost. 
These problems become even more prominent with wavelet configurations made of multiple wavelet (highpass) components, see for example the ``HHLL''-wavelet in Figure~\ref{fig:timecone}. This is because, unlike the DT-$\C$WT, by construction DWT must represent multiple ``diagonal'' directions by just one wavelet.

In figure~\ref{fig:timecone} we show the reconstruction of the same 3D ball growing over time but we now visualize the central $xz$-plane over time, resulting in a ``time-cone'' shape. Analogously to figure~\ref{fig:LLHLrecon}, we compare the reconstructions 
obtained from the coarsest scale ``LLHH''-wavelet coefficients with both dual-tree complex wavelets and (Daubechies 2) real-valued  wavelets. 
Here, real-valued wavelets (right) reveal only partially the edge perpendicular to the $zt$-plane and the region is disjoint as it shifts over time. In contrast, complex wavelets (middle), thanks also to shift-invariance, manage to represent the edge faithfully even as the singularity shifts outwards.

Notice that, in figures~\ref{fig:LLHLrecon} and~\ref{fig:timecone}, for the DT-$\C$WT wavelets from all 8 orthants are used. Furthermore, the top half of the volume is given by the 4 wavelets corresponding to orthants 1-4 (where $\I_z = +i$) and the bottom half by wavelets corresponding to orthants 5-8 ($\I_z = -i$). Hence, carefully choosing certain orthants of a particular wavelet could be used to formally analyze the geometry of the decomposed object: we leave this to future work.

\section{Applications}
\label{sec:applications}
In order to demonstrate the potentiality of the 4D DT-$\C$WT, we apply it to the inverse problem of reconstructing a volume over time, namely, 4D (3D+time) dynamic computed tomography (CT).

CT is a well known inverse problem where the inner structure of an unknown object is determined from external measurements of its X-ray attenuation intensity. This task is notoriously ill-posed, especially when only a sparse sample of measurements is available. One way to overcome ill-posedness, and therefore guarantee a stable (and unique) solution, is to add \textit{regularization} to the problem~\cite{Engl1996}. In the latest years, sparse regularization strategies, based on the paradigm that for each class of data, there exists a sparsifying representation system (such as wavelets or shearlets), have been widely used in CT applications, including dynamic CT (see~\cite{bubba2020sparse} and references therein).

Starting from the model first introduced in~\cite{bubba2020sparse} for the 2D+time case, we extend it to the 4D case, using complex wavelets rather than shearlets as a regularizer. 

\subsection{Mathematical model} 
Modern cone-beam CT scanners collect 2D projection images from given angle views. These can then be used to reconstruct a 3D volume of the interior attenuation of the targeted object. If this measurement process is then repeated over time, the object of interest can be understood as a 4D object. Given the sparse measurements and the violation of the static assumption that it is often assumed in classic CT reconstruction schemes, recovering a moving object from multiple sparse measurements over a given time period requires regularization with an appropriate representation system. Here, we use the 4D DT-$\C$WT: in analogy with the approach in~\cite{bubba2020sparse}, we are not only regularizing spatially on the 3D volume but also across time frames by considering the 3D moving volume as a 4D object.

Formally, for each time step $t = 1,...,T$, let $\f_t(x,y,z) \in \R^{N}_+$, with $N = N_xN_yN_zN_t$, be a vector representing the unknown 3D object, $\RadonD_t \in \R^{M \times N}$ a matrix modelling the tomographic cone-beam measurement process and $\m_t + \etab =: \m^{\etab}_t \in \R^M$ the data corrupted by measurement errors $\etab = \etab(t)$. To further simplify our notation we set:
\begin{equation*}
    \f = \left[ \begin{array}{c}
    \f_1 \\
    \vdots \\
    \f_T \end{array} \right], \ \RadonD = \left[ \begin{array}{ccc}
    \RadonD_1 & & \\
     & \ddots & \\
     & & \RadonD_T \end{array} \right], \ \m^{\etab} = \left[ \begin{array}{c}
    \m^{\etab}_1 \\
    \vdots \\
    \m^{\etab}_T \end{array} \right].
\end{equation*}

Then a regularized solution $\f \in \R_+^{NT}$ can be obtained by minimizing the functional
\begin{equation} \label{eq:functional}
    J(\f) = \frac{1}{2}\|\RadonD \f - \m^{\etab} \|_2^2 + \mu \|\ComplexWTD \f\|_1.
\end{equation}
Here, the regularization parameter $\mu > 0$ balances between the data mismatch over the time steps and the $\ell^1$-sparsity of the solution in the 4D dual-tree complex wavelet domain.

A robust minimization method is the primal-dual fixed point (PDFP) algorithm~\cite{chen2016primal}. Similarly to the well-known iterative soft-tresholding algorithm (ISTA), the wavelet coefficients of the iterates are soft-thresholded depending on the parameter $\mu$. Compared to ISTA, PDFP allows for additional constraints (namely the non-negativity of $\f$) and ensures convergence even when the spasifying system does not form an orthonormal basis but a frame, as it is the case with 
dual-tree complex wavelets. By using PDFP, equation~\eqref{eq:functional} can be minimized by iterating the following steps:
\begin{align} 
\begin{split}
\db^{(i+1)} &= \proj{+} \big(\f^{(i)} - \gamma (\RadonD^T \Radon \f^{(i)} - \RadonD^T \m^{\etab}) - \lambda \ComplexWTD^* \vi^{(i)} \big), \\[0.5em]
\vi^{(i+1)} &= \big(\mathbb{I} - S_{\mu \frac{\gamma}{\lambda}}\big) \big(\ComplexWTD \db^{(i+1)} + \vi^{(i)} \big), \\[0.5em]
\f^{(i+1)} &= \proj{+} \big(\f^{(i)} - \gamma (\RadonD^T \RadonD \f^{(i)} - \RadonD^T \m^{\etab}) - \lambda \ComplexWTD^* \vi^{(i+1)} \big)
\end{split}
\label{eq:PDFP}
\end{align}
$S_{\mu \frac{\gamma}{\lambda}}$ denotes the soft-thresholding operator and $\proj{+}$ is the projection onto the non-negative orthant.
The parameters $\gamma$ and $\lambda$ are bounded by properties of the functional $J$, which set a clear range for their values, while the optimal choice of $\mu$ is a notoriously difficult task. Here, we utilize an automated tuning of $\mu$ based on the a priori given desired sparsity level of the wavelet coefficients. This method was first introduced in~\cite{purisha2017controlled} using Haar wavelet regularization in traditional 2D tomography regularization and contains a detailed explanation of the automated sparsity control. Recently we applied the method to 2D+time dynamic tomography setting using shearlets~\cite{bubba2020sparse}, where we also motivated the choice of this model further.

Since the DT-$\C$WT coefficients are complex-valued it is worth noting that the soft-thresholding function in equation~\eqref{eq:PDFP} acts radially:
\[
S_{\mu'}(v) := \max \{ 0, |v|-\mu' \} e^{i\arg(v)},
\]
and component-wise when $v$ is a vector. Here, $\arg(v)$ denotes the argument of $v \in \C$.

For comparison purposes, we implemented also 4D DWT: the regularized model with 4D DWT is obtained by replacing $\ComplexWTD$ with a DWT (denoted in the following by $\WaveTD$) in equation~\eqref{eq:PDFP} and changing the values of $\lambda$ and $\mu$ accordingly. The 4D DWT is implemented by extending the 3D DWT from the Wavelet Toolbox and is available on GitHub~\cite{heikkila2021wavedec4}.

Finally, the matrices $\RadonD_t$ (and therefore $\RadonD$) simulating the geometry of cone-beam CT machine are generated using ASTRA Toolbox~\cite{van2016fast}.

\subsection{Test cases}

\begin{figure}[!tb]
\centering
\begin{tabular}{cc}
    $t=2$ & $t=4$ \\
    \includegraphics[width=.48\linewidth]{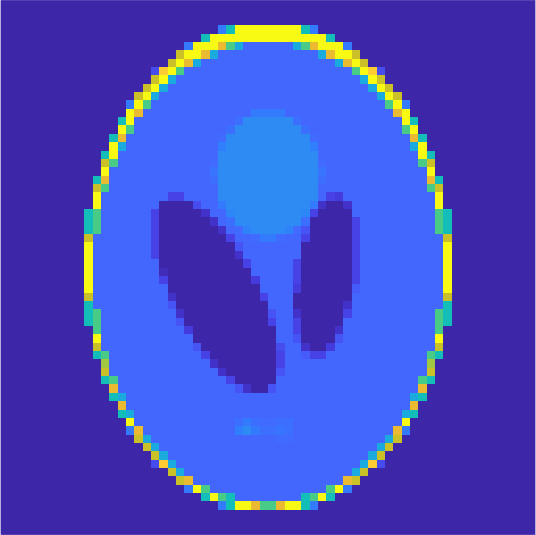} &
    \includegraphics[width=.48\linewidth]{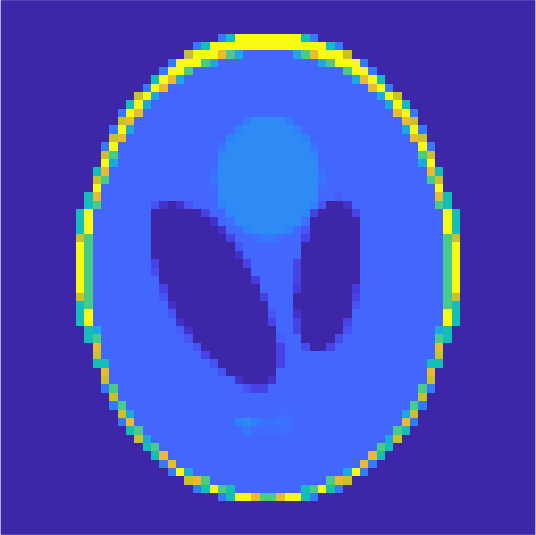} \\
    \includegraphics[width=.48\linewidth]{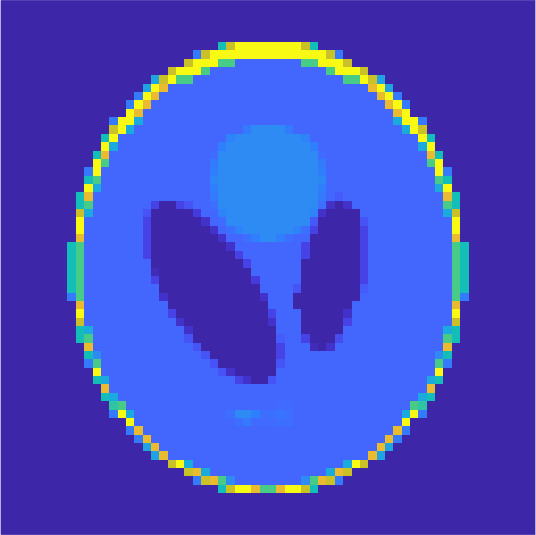} & 
    \includegraphics[width=.48\linewidth]{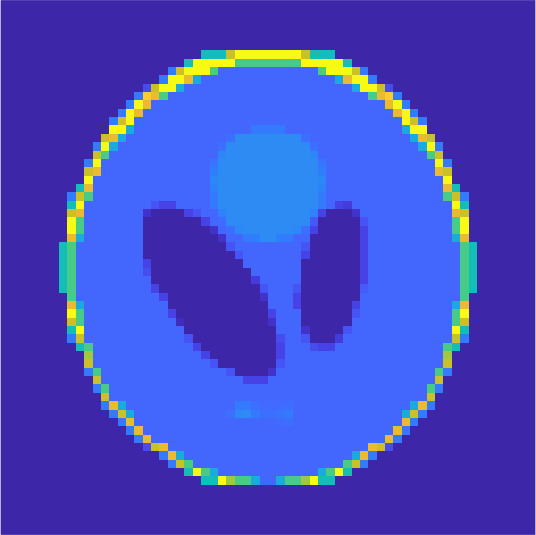} \\
    $t=6$ & $t=8$
\end{tabular}
\caption{Central horizontal slice ($z=32$) of the dynamic Shepp-Logan phantom illustrating part of the periodic deformation.}\label{fig:DSL} % Dynamic Shepp-Logan ground truth
\end{figure}

To assess the viability of 4D DT-$\C$WT regularization in sparse dynamic tomography we use two data sets which are governed by two different types of motion.
\begin{itemize}
    \item {\bf Dynamic Shepp-Logan} data is simulated by deforming a 3D version of the famous Shepp-Logan phantom~ \cite{jorgensen2010tomobox}. The deformations happen at two scales: 15 small changes evenly distributed during the simulation of each sinogram and a larger change (equivalent to 15 small changes at once) between each full measurement. This reproduces a scenario where the object is changing during a full rotation of the measurement device and there is an equally long break before the next set of measurements begins.
    The overall motion is periodic over the whole time interval and consists of simultaneous squeezing and stretching of the whole phantom in each direction. 
    
    To avoid inverse crimes the projection images are generated at twice the required resolution, down-sampled and contaminated with additive Gaussian noise with $0$ mean and $5\%$ variance. 
    Some interior slices (at $z = 32$ and different time steps $t$) of the simulated $64 \times 64 \times 64 \times 16$ object are shown in figure~\ref{fig:DSL}. The selected time steps correspond roughly to half a period of the motion.
    
    \item {\bf Gel phantom} data is from real $\mu$CT measurements of a test tube filled with agarose gel and perfused with potassium iodide contrast agent using vertical cavities in the gel body. Detailed documentation of the same setup but containing only the central slice of each projection image (for 2D + time fan-beam measurements) can be found in arXiv and the data files in Zenodo~\cite{heikkila2020gel}. Full dynamic cone-beam data used here will be made openly available in the future.
    
    The motion inside the gel phantom is only caused by the perfusing iodine and the remaining of the structure is static. However, the total intensity of the object changes at an unknown rate. To slightly increase the ill-posedness, Gaussian noise with $0$ mean and $1\%$ variance was added to the already noisy data.
\end{itemize}

To apply the automated choice for the regularization parameter, we need to determine the a priori level of sparsity~\cite{purisha2017controlled}. 
The desired sparsity level for the dynamic Shepp-Logan data was calculated from the known 4D object and was chosen to be $d_{\ComplexWTD} = 0.6$ for the DT-$\C$WT and $d_{\WaveTD} = 0.5$ for the DWT. For the gel phantom data we used as ``ground truth'' a high quality reconstruction obtained with the FDK-algorithm~\cite{feldkamp1984practical} using 360 projection angles and no additional noise. The desired sparsity levels were chosen to be $d_{\ComplexWTD} = 0.6$ for the DT-$\C$WT and $d_{\WaveTD} = 0.6$ for the DWT. 

Notice that these were also the reference objects used for the numerical error estimates reported in table~\ref{tab:errors}.

\subsection{Results} \label{ssec:results}

\begin{figure*}[!tb] % test 4.1
\centering
\setlength{\tabcolsep}{0.1em}
\begin{tabular}{rccccl}
    & $t=2$ & $t=4$ & $t=6$ & $t=8$ & \\
    \rotatebox[origin=l]{90}{\hspace{0.035\linewidth} $xy$-plane, DWT} &
    \includegraphics[width=.23\linewidth]{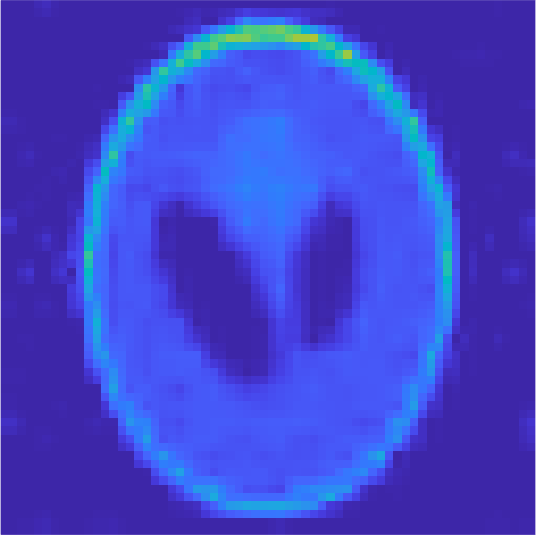} &
    \includegraphics[width=.23\linewidth]{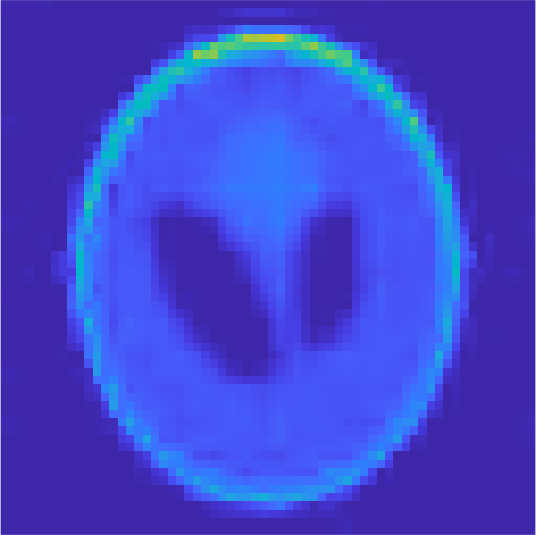} & 
    \includegraphics[width=.23\linewidth]{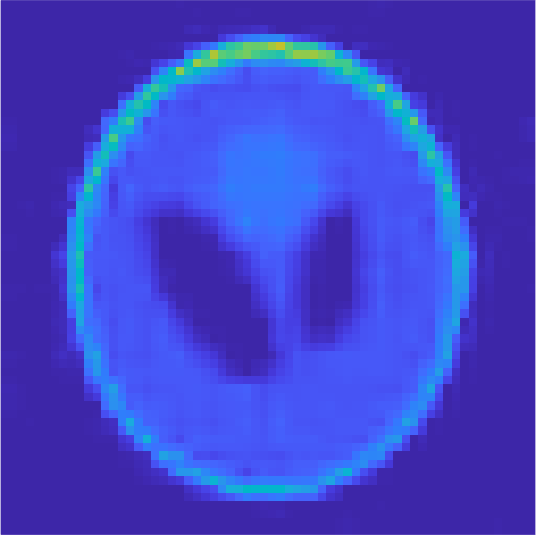} & 
    \includegraphics[width=.23\linewidth]{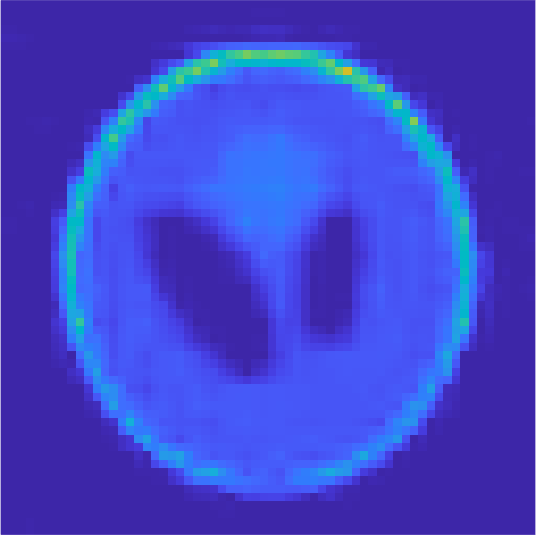} & 
    \includegraphics[height=.23\linewidth]{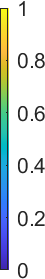} \\ % Colorbar
    \rotatebox[origin=l]{90}{\hspace{0.03\linewidth} $xy$-plane, DT-$\C$WT} &
    \includegraphics[width=.23\linewidth]{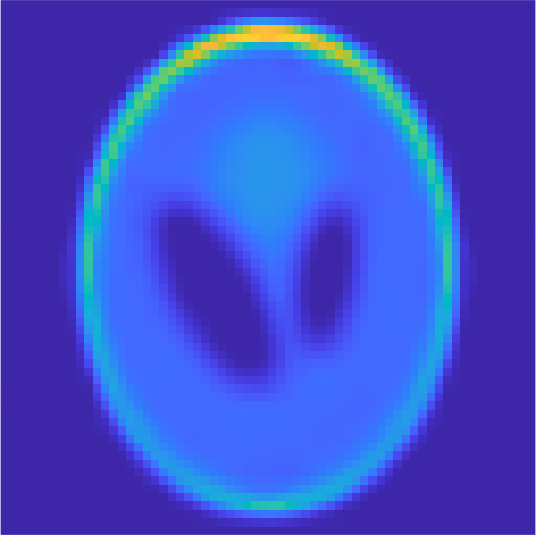} &
    \includegraphics[width=.23\linewidth]{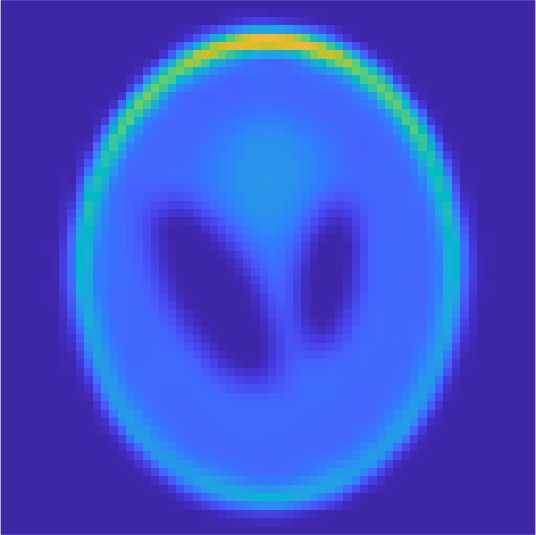} & 
    \includegraphics[width=.23\linewidth]{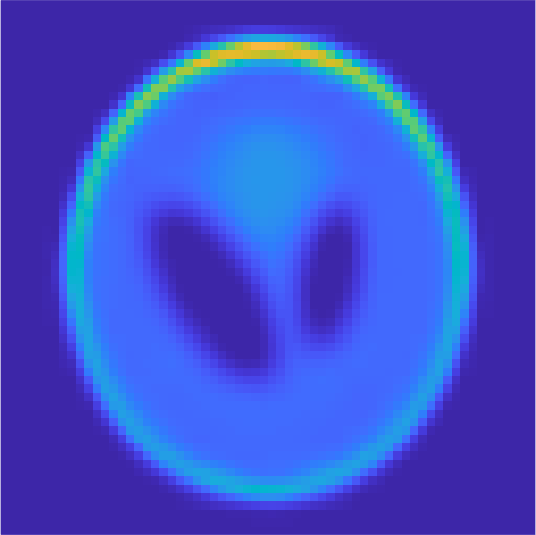} & 
    \includegraphics[width=.23\linewidth]{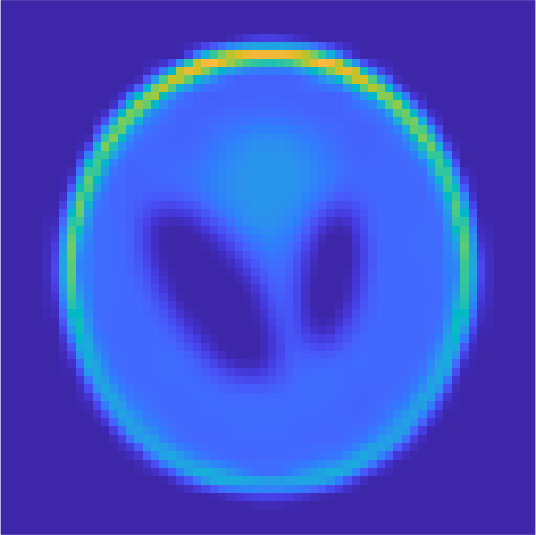} & 
    \includegraphics[height=.23\linewidth]{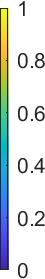} \\ % Colorbar
\end{tabular}
\caption{Horizontal slices at $z=32$ of the dynamic Shepp-Logan phantom reconstruction at various time steps and using 30 projections. Regularization with Daubechies 2 real-valued wavelets (DWT, top row) and dual-tree complex wavelets (DT-$\C$WT, bottom row).}\label{fig:DSL30ang} % Dynamic Shepp-Logan
\end{figure*}

Reconstructions from the dynamic Shepp-Logan data can be seen in figure~\ref{fig:DSL30ang} (using DWT and DT-$\C$WT). Similarly to figure~\ref{fig:DSL}, we show the horizontal ($xy$-plane) slice at height $z = 32$ of selected time steps.

\begin{figure*}[!tb] % test 4
\centering
\setlength{\tabcolsep}{0.1em}
\begin{tabular}{rccccl}
    & $t=4$ & $t=8$ & $t=12$ & $t=16$ & \multirow{3}{*}{
    \includegraphics[height=.42\linewidth]{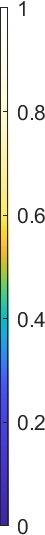}} \\ % Multirow colorbar
    \rotatebox[origin=l]{90}{\hspace{0.06\linewidth} $xy$-plane} &
    \includegraphics[width=.23\linewidth]{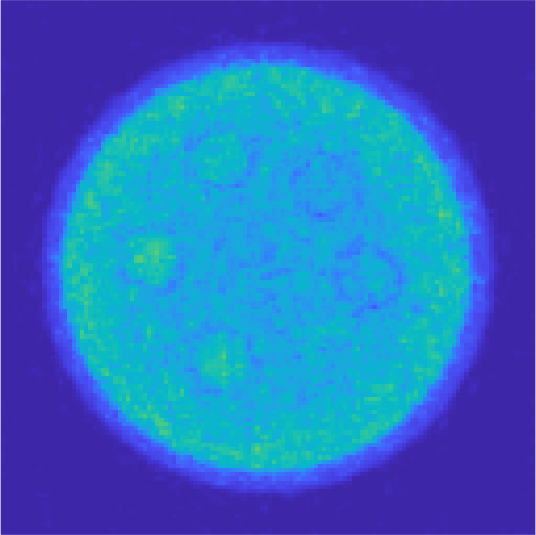} &
    \includegraphics[width=.23\linewidth]{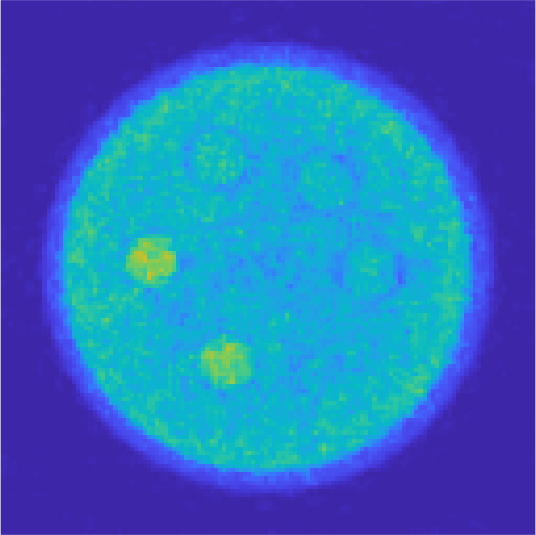} & 
    \includegraphics[width=.23\linewidth]{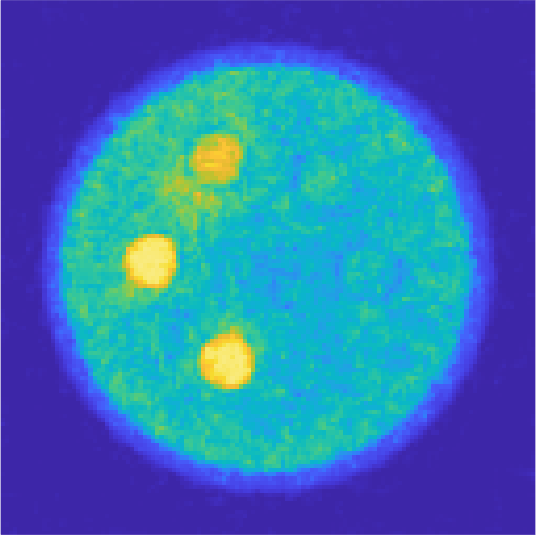} & 
    \includegraphics[width=.23\linewidth]{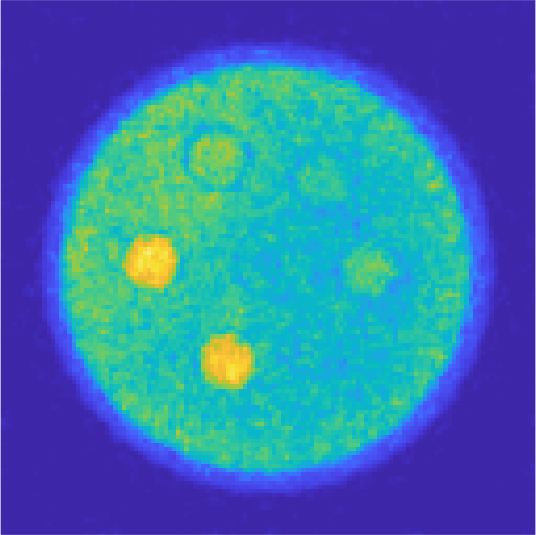} & \\
    \rotatebox[origin=l]{90}{\hspace{0.04\linewidth} $xz$-plane} &
    \includegraphics[width=.23\linewidth]{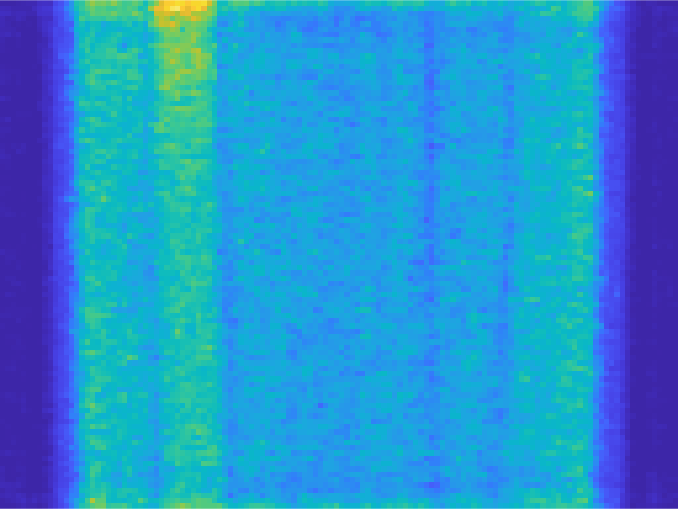} &
    \includegraphics[width=.23\linewidth]{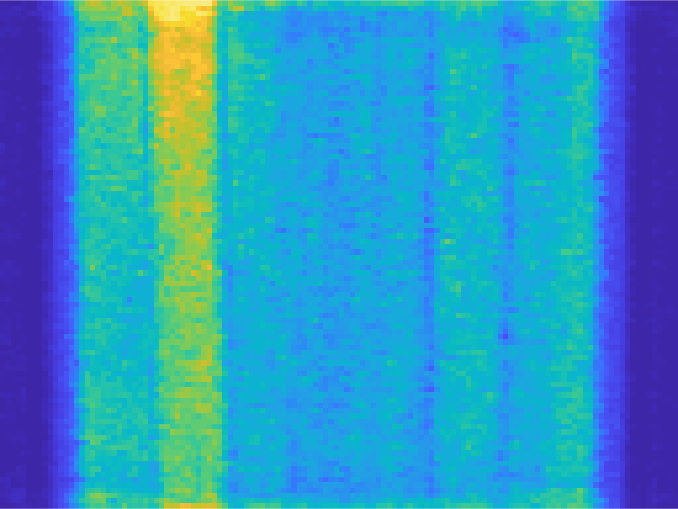} &
    \includegraphics[width=.23\linewidth]{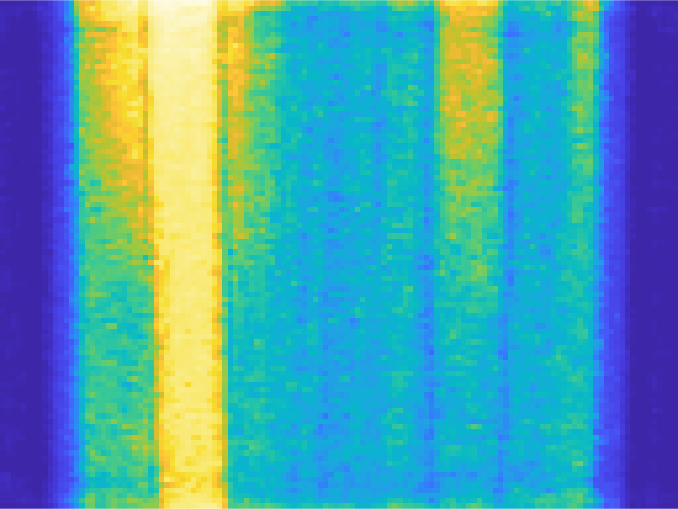} & 
    \includegraphics[width=.23\linewidth]{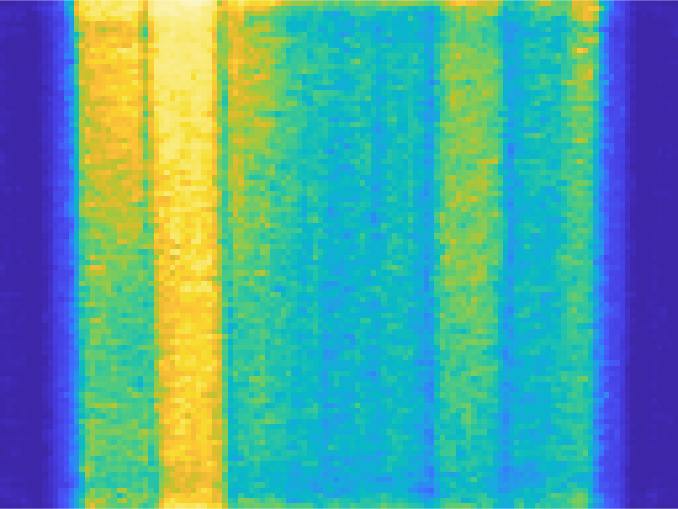} &
\end{tabular}
\caption{Horizontal slices at $z=64$ (top row) and vertical slices at $y=64$ (bottom row) of gel phantom reconstruction at various time steps, using 30 projections. Regularization with real-valued Daubechies 2 wavelets (DWT).}\label{fig:5pt30angW} % 5 point test (gel phantom)
\end{figure*}

\begin{figure*}[!tb] % test 4
\centering
\setlength{\tabcolsep}{0.1em}
\begin{tabular}{rccccl}
    & $t=4$ & $t=8$ & $t=12$ & $t=16$ & \multirow{3}{*}{
    \includegraphics[height=.42\linewidth]{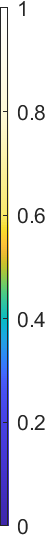}} \\ % Multirow colorbar
    \rotatebox[origin=l]{90}{\hspace{0.06\linewidth} $xy$-plane} &
    \includegraphics[width=.23\linewidth]{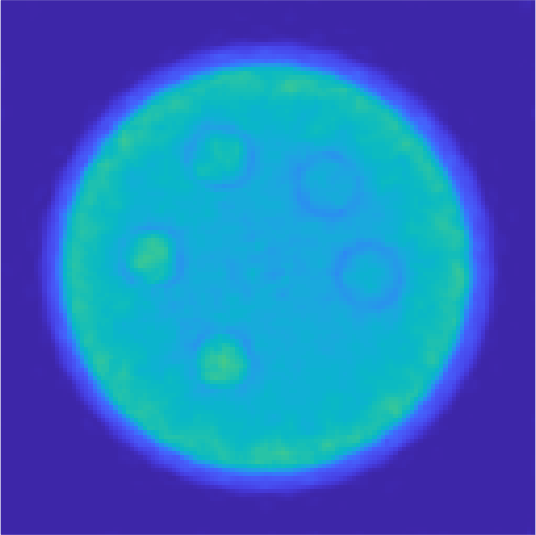} &
    \includegraphics[width=.23\linewidth]{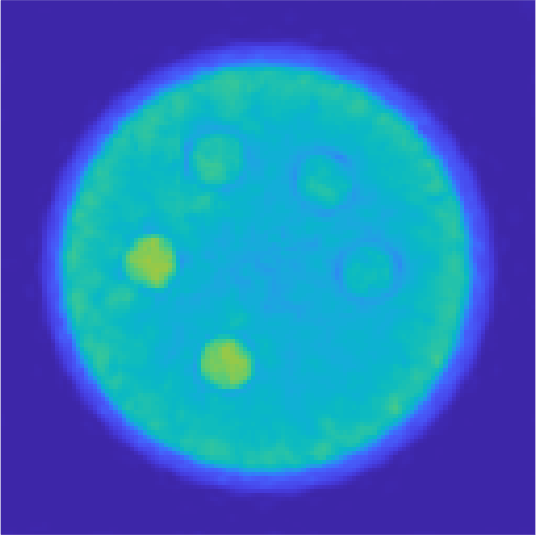} & 
    \includegraphics[width=.23\linewidth]{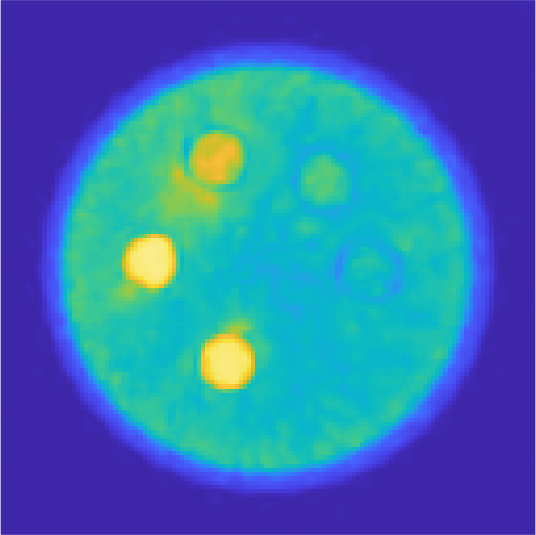} & 
    \includegraphics[width=.23\linewidth]{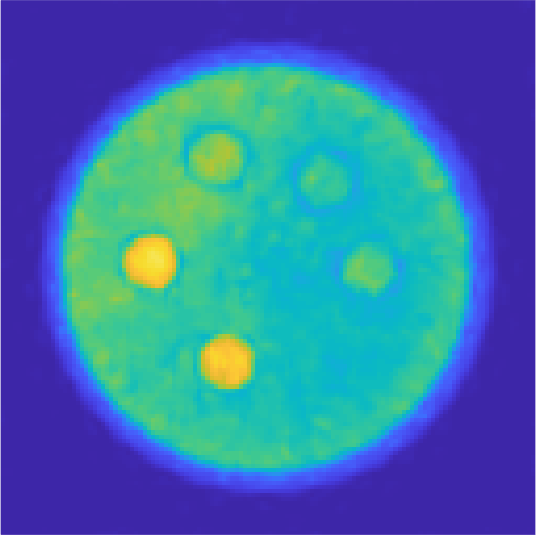} & \\
    \rotatebox[origin=l]{90}{\hspace{0.04\linewidth} $xz$-plane} &
    \includegraphics[width=.23\linewidth]{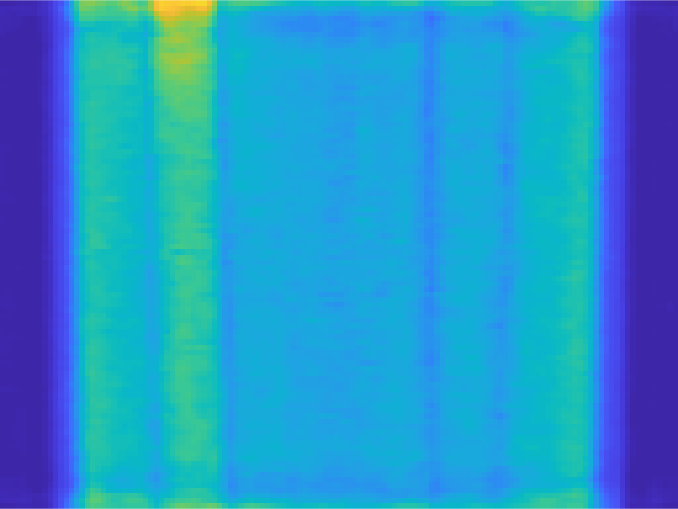} &
    \includegraphics[width=.23\linewidth]{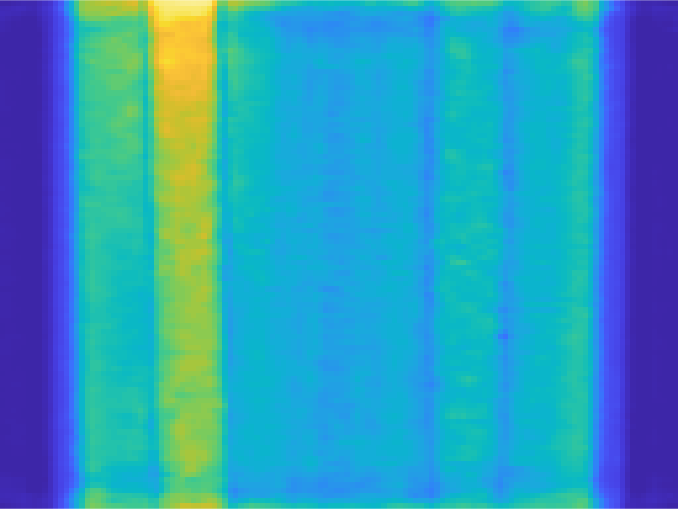} &
    \includegraphics[width=.23\linewidth]{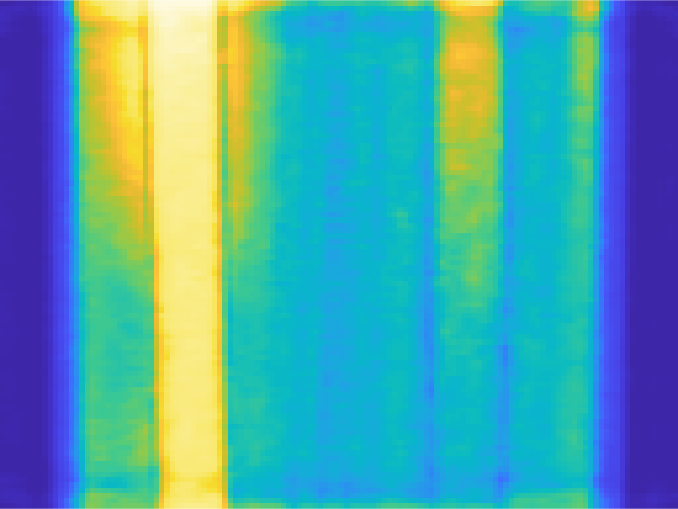} & 
    \includegraphics[width=.23\linewidth]{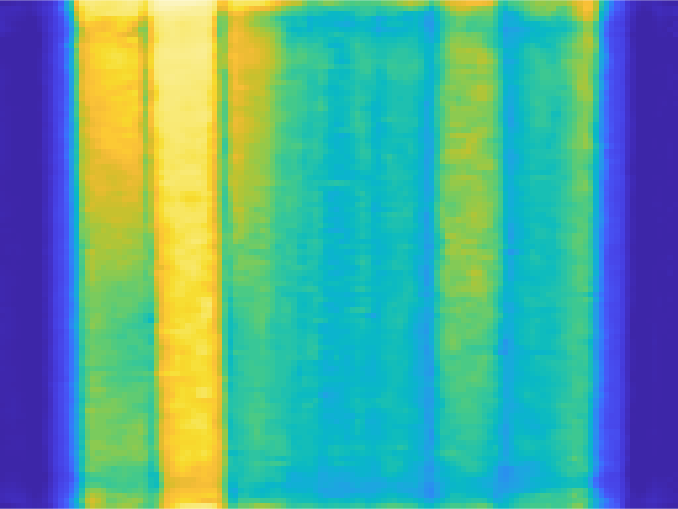} &
\end{tabular}
\caption{Horizontal slices at $z=64$ (top row) and vertical slices at $y=64$ (bottom row) of gel phantom reconstruction at various time steps, using 30 projections. Regularization with dual-tree complex wavelets (DT-$\C$WT).}\label{fig:5pt30angCWT} % 5 point test (gel phantom)
\end{figure*}

Reconstructions from the gel phantom data are reported in figures \ref{fig:5pt30angW} (using DWT) and \ref{fig:5pt30angCWT} (using DT-$\C$WT). In each column we show 2D slices from selected time steps: on the top row there is the horizontal ($xy$-plane) slice at height $z = 64$ and on the bottom row there is the vertical ($xz$-plane) slice at $y=64$. Both reconstructions use 30 projection angles and were originally of size $128 \times 128 \times 128 \times 16$ but have been cropped vertically (along $z$-axis) to size $128 \times 128 \times 96 \times 16$. This is done to avoid artifacts caused by the phantom extending vertically outside the measured X-ray cone.

In addition to the visual comparisons some numerical error estimates for both data are provided. Relative $\ell^2$-norm error and peak signal-to-noise (PSNR) ratios of the 4D reconstruction are listen in table~\ref{tab:errors}. We also consider the Haar-wavelet perceptual similarity index (HPSI)~\cite{reisenhofer2018} of the horizontal slices of the dynamic Shepp-Logan (seen in figure~\ref{fig:DSL30ang}) and the vertical slices (bottom row in figures~\ref{fig:5pt30angW} and~\ref{fig:5pt30angCWT}) of the gel phantom. We then calculate the mean value over all 16 time steps. For all numerical error estimates the gel phantom reconstructions are cropped vertically to more fairly evaluate the regularization without the cone-beam geometry artifacts.

Comparing the quality of reconstructions in figures \ref{fig:DSL30ang}, \ref{fig:5pt30angW} and \ref{fig:5pt30angCWT} we notice that overall the dual-tree complex wavelets perform better at preserving details whilst also denoising the reconstructions. For example, in figure~\ref{fig:DSL30ang} it can be seen that reconstructions with Dauchechies 2 wavelet regularization (top row) are clearly noisier and the outer boundary is not nearly as well preserved as with DT-$\C$WT regularization (especially at $t = 8$). This can be taken as evidence that complex-valued wavelets are better at preserving these features thanks to their shift-invariance and directional sensitivity.

The differences with the gel phantom reconstructions are not as remarkable but, again, the DT-$\C$WT regularized solution has noticeably less noise and especially the edges of the vertical cavities (dark blue circles in the $xy$-plane images) are better preserved. The bright iodine (in yellow) is well reconstructed by both methods.

The numerical error estimates in table~\ref{tab:errors} provide less insights but seem to favour the DT-$\C$WT reconstructions with the exception of the mean HPSI of the dynamic Shepp-Logan data where DWT obtains slightly better values. Notice also that based on table~\ref{tab:timetable} the computational cost of the DT-$\C$WT regularization is ``only'' about 10-times larger than the DWT regularization compared to the roughly 16-fold increase in computations of the wavelet transform itself.

Finally, we incidentally mention that the inclusion of the third spatial direction seems to improve the quality of robustness of the reconstructions compared to the similar 2D + time setup in~\cite{bubba2020sparse}. While the angular sampling is definitely sparse (just 30 projections), this does not affect the $z$-direction which provides additional robustness and seems to decrease to some extent the ill-posedness of the problem.

\begin{table}
    \centering
    \caption{Numerical error estimates of the different reconstructions.}
    \label{tab:errors}
    \begin{tabular}{ll m{4em} l m{4em}}
        & & Relative error & PSNR & Mean HPSI \\
       \multirow{2}{2cm}{Dynamic Shepp-Logan} & DT-$\C$WT & \textbf{40.3\%} & \textbf{22.66} & 0.603 \\
        & DWT & 44.8\% & 21.73 & \textbf{0.643} \\ \hline 
        \multirow{2}{2cm}{Gel phantom} & DT-$\C$WT & \textbf{9.30\%} & \textbf{30.86} & \textbf{0.637} \\
        & DWT & 10.83\% & 29.54 & 0.621 \\ 
    \end{tabular}
\end{table}

\begin{table}
    \centering
    \caption{Number of iterations and computational times of the different reconstructions.}
    \label{tab:timetable}
    \begin{tabular}{ll c c c}
        & & \multirow{2}{*}{Iterations} & \multicolumn{2}{c}{Time ($s$)} \\
        & & & total & per iter. \\
        \multirow{2}{2cm}{Dynamic Shepp-Logan} & DT-$\C$WT & $70$ & $1531$ & $21.9$ \\
        & DWT & $63$ & $137$ & $2.2$ \\ \hline
        \multirow{2}{2cm}{Gel phantom} & DT-$\C$WT & $61$ & $11357$ & $186.2$ \\
        & DWT & $54$ & $902$ & $16.7$ \\
    \end{tabular}
\end{table}

\section{Conclusions}
\label{sec:conclusions}

In this paper we introduced the 4D DT-$\C$WT and explored its use to address the inverse problem of reconstructing a volume evolving over time from dynamic tomographic data. Our analysis speaks in favor of this type of representation to address space-time problems thanks to its simple implementations and  strong theoretical properties. Our results show a potential for 4D complex wavelets to be competitive in a numerical framework even when compared to other (more) refined multidimensional systems.

\printbibliography
%\bibliography{IEEEabrv,references}

\end{document}